\documentclass[10pt]{amsart}

\usepackage{amsmath,amssymb}
\usepackage{ulem}
\usepackage[dvipdfmx]{graphicx,xcolor} 
\usepackage{amscd} 
\usepackage{array} 
\usepackage{float} 

\newtheorem{theorem}{Theorem}[section]
\newtheorem{lemma}[theorem]{Lemma}
\newtheorem{proposition}[theorem]{Proposition}
\newtheorem{corollary}[theorem]{Corollary}

\theoremstyle{definition}

\theoremstyle{remark}
\newtheorem{remark}[theorem]{Remark}

\numberwithin{equation}{section}

\makeatletter
\@namedef{subjclassname@2020}{\textup{2020} Mathematics Subject Classification}
\makeatother

\begin{document}

\title[Attraction rates]
{Attraction rates for iterates of a superattracting skew product}  

\author[K. Ueno]{Kohei Ueno}
\address{Daido University, Nagoya 457-8530, Japan}
\curraddr{}
\email{k-ueno@daido-it.ac.jp}
\thanks{This work was supported by the Research Institute for Mathematical Sciences, 
           an International Joint Usage/Research Center located in Kyoto University.}



\date{}
\dedicatory{}

\maketitle

\begin{abstract}
Let $f(z,w)=(p(z),q(z,w))$ be a holomorphic skew product
with a superattracting fixed point at the origin.
In the previous paper
we have succeeded to specify a dominant term of $q$ 
by the order of $p$ and the Newton polygon of $q$ 
and to construct a B\"{o}ttcher coordinate
on an invariant wedge.
By using the same idea and terminologies,  
we give inequalities on attraction rates
for the vertical dynamics of $f$
in this paper.
The results hold not only for the superattracting case,
but for all the other cases.
\end{abstract}

\section{Introduction}

Let $f : (\mathbb{C}^2, 0) \to (\mathbb{C}^2, 0)$ be a holomorphic germ
with a superattracting fixed point at the origin
and let $f^n$ be the $n$-th iterate of $f$.
We define $c(f)$ as the smallest degree of any term 
in the Taylor expansion of $f$ in local coordinates,
which is independent of the choice of coordinates,
and we call $c(f)$ the attraction rate of $f$.
From the viewpoint of complex dynamics, 
it is important to study the behavior of the attraction rates for iterates of $f$,
because it gives a measure of the rate 
at which nearby points are attracted to the origin under iteration.
Moreover,
the limit $c_{\infty} = \lim_{n \to \infty} \sqrt[n]{c(f^n)}$,
which we call the asymptotic attraction rate of $f$,
measures the growth of the sequence $\{ c(f^n) \}_{n \geq 1}$ of attraction rates.
Favre and Jonsson \cite{fj} proved that 
$c_{\infty}$ is a quadratic integer and
there exists $D \in ( 0,1 ]$ such that 
$D c_{\infty}^n \leq c(f^n) \leq c_{\infty}^n$ for any $n \geq 1$.
This result is derived from their result on normal forms of $f$,
which is obtained by blow-ups 
and applied to construct a pluriharmonic function 
with the adequate invariance property. 
By using valuative techniques similar to those developed by Favre and Jonsson,
Gignac and Ruggiero \cite{gr} proved that 
the sequence $\{ c(f^n) \}_{n \geq 1}$ eventually satisfies 
an integral linear recursion relation,
which, up to replacing $f$ by an iterate,
can be taken to have order at most two. 

In this paper 
we are concerned with skew products, 
and provide equalities and inequalities on the attraction rates for the vertical dynamics,
using the same idea and terminologies as in the previous paper \cite{u}.
A holomorphic germ of the form $f(z,w)=(p(z),q(z,w))$ is called a skew product.
We assume that it has a fixed point at the origin;
thus we have the Taylor expansions
\[
p(z) = a_{\delta} z^{\delta} + O(z^{\delta + 1})
\text{ and }
q(z,w) = \sum_{i+j \geq 1} b_{ij} z^{i} w^{j},
\]
where $a_{\delta} \neq 0$ and $\delta \geq 1$,
in local coordinates.
We define $c(f)$, $c(p)$ and $c(q)$ as 
the smallest degrees of any term in this Taylor expansion of $f$, $p$ and $q$,
respectively.
Then 
$c(f) = \min \{ c(p), c(q) \}$, 
$c(p) = \delta$ and 
\[
c(q) = \min \{ i+j \, | \, b_{ij} \neq 0 \}.
\]
Let $f^n(z,w) = (p^n(z), Q^n(z,w))$.
Then $c(f^n) = \min \{ c(p^n), c(Q^n) \}$ and $c(p^n) = \delta^n$.
Our aim is to provide an estimate on the attraction rate $c(Q^n)$,  
which is deeply related to the Newton polygon of $Q^n$
as stated later.
We remark that,
although an estimate on $c(Q^n)$ implies that on $c(f^n)$,
the opposite direction is not true in general.


The fundamental properties of the dynamics of polynomial skew products is
well studied and summarized in \cite{j} and \cite{fg}.
Lilov \cite{l} studied the local and semi-local dynamics of 
holomorphic skew products near a superattracting invariant fiber.
See also \cite{ps} and \cite{pr}
for the dynamics of skew products near an invariant fiber
of different types.
Whereas the dynamics of skew products is mild 
in the sense that it can be regarded as 
the intermediate between the one and two dimensional dynamics,
it has very complicated aspects and
exhibits new phenomenon of the dynamics in dimension two.
For example,
Astorg {\it et al} \cite{abdpr} 
found polynomial skew products with wandering domains,
which is the first example of non-invertible
polynomial maps with wandering domains. 
See also \cite{abp} and \cite{ab}.
Furthermore,
polynomial skew products are used in \cite{d} and \cite{t}
to construct robust bifurcations.

Let us first recall our previous result in \cite{u}. 
Assuming that the origin is superattracting,
we have succeeded  
to show the existence of a dominant term $b_{\gamma d} z^{\gamma} w^d$ of $q$ and
to construct a B\"{o}ttcher coordinate that conjugates $f$ to 
the monomial map $f_0$ on the wedge $U$,
where $f_0(z,w) = (a_{\delta} z^{\delta}, b_{\gamma d} z^{\gamma} w^d)$
and 
$U = \{ |z|^{l_1 + l_2} < r^{l_2} |w|, |w| < r|z|^{l_1} \}$
for some rational numbers $0 \leq l_1 < \infty$ and $0 < l_2 \leq \infty$
and for small $r > 0$.
The bidegree $(\gamma, d)$ and 
the rational numbers $l_1$ and $l_2$ are determined by
the order of $p$ and the Newton polygon of $q$.
We define the Newton polygon $N(q)$ of $q$ as 
the convex hull of the union of $D(i,j)$ with $b_{ij} \neq 0$,
where $D(i,j) = \{ (x,y) \, | \, x \geq i, y \geq j \}$. 
Let $(n_1, m_1)$, $(n_2, m_2)$, $\cdots, (n_s, m_s)$ be the vertices of $N(q)$,
where $n_1 < n_2 < \cdots < n_s$ and $m_1 > m_2 > \cdots > m_s$.
Let $T_k$
be the $y$-intercept of the line $L_k$ passing through 
the vertices $(n_k, m_k)$ and $(n_{k+1}, m_{k+1})$
for each $1 \leq k \leq s-1$.

\vspace{2mm}
\begin{itemize}
\item[Case 1] 
If $s = 1$, then $N(q)$ has a unique vertex, which is denoted by $(\gamma, d)$, \\
and we define
$l_1 =  l_2^{-1} = 0$.
Hence $U = \{ |z| < r, |w| < r \}$. 
\end{itemize}
\vspace{1mm}

If $s = 1$, then 
$b_{\gamma d} z^{\gamma} w^{d}$ is clearly the dominant term of $q$ 
and the result is classical. 
Difficulties appear when $s > 1$,
which is divided into the following three cases. 

\vspace{2mm}
\begin{itemize}
\item[Case 2] 
  If $s > 1$ and $\delta \leq T_{s-1}$, then  we define
  \[
  (\gamma, d) = (n_s, m_s), \ 
  l_1 = \frac{n_s - n_{s-1}}{m_{s-1} - m_s} 
  \text{ and } l_2^{-1} = 0.
  \]
  Hence $U = \{ |z| < r, |w| < r |z|^{l_1} \}$.
\item[Case 3] 
  If $s > 1$ and $T_1 \leq \delta$, then  we define
  \[
  (\gamma, d) = (n_1, m_1), \
  l_1 = 0 
  \text{ and } l_2 = \frac{n_2 - n_1}{m_1 - m_2}.
  \]
  Hence $U = \{ |z|^{l_2} < r^{l_2} |w|, |w| < r \} = \{ r^{-l_2} |z|^{l_2} < |w| < r \}$.  
\item[Case 4]
  If $s > 2$ and $T_k \leq \delta \leq T_{k-1}$ for some $2 \leq k \leq s-1$, then  we define
  \[
  (\gamma, d) = (n_k, m_k), \ 
  l_1 = \frac{n_k - n_{k-1}}{m_{k-1} - m_k} 
  \text{ and } l_1 + l_2= \frac{n_{k+1} - n_k}{m_k - m_{k+1}}.
  \] 
  Hence  $U = \{ r^{-l_2} |z|^{l_1 + l_2} < |w| < r|z|^{l_1} \}$.
\end{itemize}
\vspace{1mm}

Note that $\gamma > 0$ for Case 2,
$d > 0$ for Case 3, and
$\gamma d > 0$ for Case 4 by the setting.
The rational numbers
$- l_1^{-1}$ and $- (l_1 + l_2)^{-1}$ are the slopes of the lines $L_{k-1}$ and $L_k$
for Case 4 and
the same correspondence holds for all the cases
if we define $L_0 = \{ x = n_1 \}$ and $L_s = \{ y = m_s \}$.

A precise statement of our previous result is the following.

\begin{theorem}[Lemmas 1.1 and 1.7 and Theorems 1.2 and 1.8 in \cite{u}] \label{Bottcher coordinate}
If $\delta \geq 2$ and $d \geq 2$ 
or if $\delta \geq 2$, $d=1$ and $\delta \neq T_{k}$ for any $k$, 
then
$f$ preserves $U$ and 
there is a biholomorphic map $\phi$ defined on $U$  
that conjugates $f$ to $f_0$
for small $r > 0$.
\end{theorem}

The map $\phi$ is called the B\"{o}ttcher coordinate for $f$ on $U$
and constructed as the limit of 
the compositions of $f_0^{-n}$ and $f^n$. 
As detailed versions using intervals of weights,  
we also exhibit Theorems \ref{Bottcher for case 2}, 
\ref{Bottcher for case 3} and 
\ref{Bottcher for case 4}  
for Cases 2, 3 and 4,
respectively, in Section 2.

Let us next state our results on the attraction rates.
We have found that the same idea as in \cite{u} can be applied to
study the attraction rate of $Q^n$,
not only for the superattracting case but for all the other cases. 
For Case 1, it is clear that
the dominant term of $Q^n$ is $z^{\gamma_n} w^{d^n}$ 
and so $c(Q^n) = \gamma_n + d^n$ for any $n \geq 1$,
where $\gamma_n = \gamma (\delta^{n-1} + \delta^{n-2} d + \cdots + d^{n-1})$
and $\gamma_{1} = \gamma$.
Hence $c(f^n) = \min \{ \delta^n, \gamma_n + d^n \}$
and so $c_{\infty} = \min \{ \delta, d \}$.
Difficulties appear for Cases 2, 3 and 4,
and the situation differs whether $d > 0$ or $d = 0$ for Case 2.
To overcome the difficulties,
it is useful to consider the following quantity:
\[
w(q) = w_l (q)= \min \{ i+lj \, | \, b_{ij} \neq 0 \}.
\]  
In fact, 
we obtain the following equalities on $w(Q^n)$.

\begin{theorem} \label{main thm of weights}
It follows for any $n \geq 1$ that $Q^n$ contains the term $z^{\gamma_n} w^{d^n}$ and
\begin{enumerate}
\item[$(1)$]  $w_{l_1} (Q^n) = \gamma_n + l_1 d^n$ if $d>0$ for Case 2, 
\item[$(2)$]  $w_{l_2} (Q^n) = \gamma_n + l_2 d^n$ for Case 3, and
\item[$(3)$]  $w_{l_1} (Q^n) = \gamma_n + l_1 d^n$ and 
        $w_{l_1 + l_2} (Q^n) = \gamma_n + (l_1 + l_2) d^n$ for Case 4.
\end{enumerate}
Moreover, 
it follows for any $n \geq 1$ that
\begin{enumerate}
\item[$(4)$] $w_{l_1} (Q^n) = \gamma_n$ if $d=0$ for Case 2.
\end{enumerate}
\end{theorem}

Although we omit the coefficient of the term $z^{\gamma_n} w^{d^n}$ 
in the statement for simplicity,
it is equal to
$a_{\delta}^{\gamma_{n-1} + \gamma_{n-2} + \cdots + \gamma_{2} + \gamma_{1}}
b_{\gamma d}^{d^{n-1} + d^{n-2} + \cdots + d + 1}$.
Whereas $Q^n$ contains the term $z^{\gamma_n}$ for any $n \geq 1$
if $d = 0$ and $\delta < T_{s-1}$,
the term $z^{\gamma_n}$ may vanish
if $d = 0$ and $\delta = T_{s-1}$ for Case 2.
However,
the equality (4) in the theorem 
follows from the existence of the other vertex of $N(Q^n)$
that should be previous to $(\gamma_n, 0)$
if the term $z^{\gamma_n}$ did not vanish
for any $n \geq 1$.
As detailed versions using intervals of weights, 
we also provide 
Theorems \ref{main thm of weights when d>0},
\ref{main thm of weights when d=0},
\ref{main thm of weights for case 3} and 
\ref{main thm of weights for case 4} for each cases.

Because $w(Q^n)$ coincides with the minimum $x$-intercept of 
the lines with slope $- l^{-1}$ 
that intersect the Newton polygon $N(Q^n)$,
Theorem \ref{main thm of weights} 
shows us the shape of $N(Q^n)$.
More precisely,
$N(Q^n)$ is included in the upper-right region
that is surrounded by the two lines with slopes $-l_1^{-1}$ and $-(l_1+l_2)^{-1}$,
which intersect at $(\gamma_n, d^n)$.
Because $c(Q^n)$ coincides with the minimum $x$-intercept of 
the lines with slope $-1$ that intersect $N(Q^n)$,
Theorem \ref{main thm of weights} induces 
the following inequalities on the attraction rate $c(Q^n)$.

\begin{theorem} \label{main thm on attr rates}
It follows for any $n \geq 1$ that $c(Q^n) \leq \gamma_n + d^n$ and
\begin{enumerate}
\item[$(1)$] $\min \{ l_1^{-1}, 1 \} \gamma_n + d^n \leq c(Q^n)$
if $d>0$ for Case 2,
\item[$(2)$] $\gamma_n + \min \{ l_2, 1 \} d^n \leq c(Q^n)$
for Case 3, and
\item[$(3)$] $\min \{ l_1^{-1}, 1 \} \gamma_n + \min \{ l_1 + l_2, 1 \} d^n \leq c(Q^n)$
for Case 4.
\end{enumerate}
On the other hand,
it follows for any $n \geq 1$ that
\begin{enumerate}
\item[$(4)$] $\min \{ l_1^{-1}, 1 \} \gamma_n \leq c(Q^n) \leq \max \{ l_1^{-1}, 1 \} \gamma_n$
if $d=0$ for Case 2.
\end{enumerate}
\end{theorem}

This theorem is restated as Corollaries \ref{main cor on attr rates when d>0},
\ref{main cor on attr rates when d=0},
\ref{main cor on attr rates for case 3} and 
\ref{main cor on attr rates for case 4} for each cases. 
Moreover,
we can improve these inequalities
by investigating the vertices of $N(Q^n)$
that are previous and/or next to $(\gamma_n, d^n)$;
see Theorems \ref{main thm of attr rates when d>0},
\ref{main thm of attr rates when d=0 and delta < T},
\ref{main thm of attr rates when d=0 and delta = T},
\ref{main thm of attr rates for case 3} and 
\ref{main thm of attr rates for case 4} for improved versions.

Let $\alpha = \gamma/(\delta - d)$ when $\delta \neq d$.
As a corollary of Theorem \ref{main thm on attr rates},
we obtain the following inequalities on the attraction rate $c(f^n)$.

\begin{corollary} \label{main cor on attr rates}
Let $\gamma d > 0$.
Then $c_{\infty} = \delta$ and it follows for any $n \geq 1$ that
\begin{enumerate}
\item[$(1)$] $\alpha \delta^n \leq c(f^n) < \delta^n$ 
if $\delta > d$ and $\alpha < 1$, and 
\item[$(2)$] $c(f^n) = \delta^n$ 
if $\delta > d$ and $\alpha \geq 1$ or if $\delta \leq d$
\end{enumerate}
for Cases 2, 3 and 4.
On the other hand,
it follows for any $n \geq 1$ that
\begin{enumerate}
\item[$(3)$] $c_{\infty} = \delta$ and
$\min \{ 1, \gamma/\delta, l_1^{-1} \gamma/\delta \} \delta^n \leq c(f^n) \leq \delta^n$
if $d = 0$ for Case 2, and
\item[$(4)$] $c_{\infty} = d$ and 
$\min \{ 1, l_2 \} d^n \leq c(f^n) \leq d^n$
if $\gamma = 0$ for Case 3.
\end{enumerate}
\end{corollary}

Using the same idea and terminologies as in \cite{u},
we constructed a pluriharmonic function in \cite{u - psh funs}
that describes the vertical dynamics well
for a superattracting skew product.
Corollary \ref{main cor on attr rates} 
is used to show that 
our pluriharmonic function is more useful to study the vertical dynamics 
than the pluriharmonic function constructed in \cite{fj}.

The organization of this paper is as follows.
In Section 2
we review the related definitions and results in \cite{u}.
More precisely,
we recall the definitions of the intervals of weights
and illustrate 
Theorems \ref{Bottcher for case 2},
\ref{Bottcher for case 3} and 
\ref{Bottcher for case 4},
detailed versions of 
Theorem \ref{Bottcher coordinate},
by blow-ups. 
Although we do not use blow-ups to prove the theorems,
they are useful to explain the theorems
when the weights are integers.
We then prove 
Theorems \ref{main thm of weights} and \ref{main thm on attr rates}
or, more precisely,
detailed versions of 
Theorem \ref{main thm of weights} and
improved versions of
Theorem \ref{main thm on attr rates},
for Case 2 when $d>0$,
Case 2 when $d=0$, 
Case 3 and Case 4
in Sections 3, 4, 5 and 6,
respectively.
Finally,
we deduce 
Corollary \ref{main cor on attr rates}
from Theorem \ref{main thm on attr rates}
in Section 7.

\section{Intervals of weights and Blow-ups} 

In this section
we give a summary of our previous results in \cite{u}:
we introduce intervals of weights,
explain benefits of the intervals in terms of blow-ups, 
and state the results on B\"{o}ttcher coordinates in terms of the intervals. 
We deal with Cases 2, 3 and 4 in Sections 2.1, 2.2 and 2.3,
respectively.
The intervals for Cases 2 and 3 are also used to 
describe equalities on $w(Q^n)$ in Sections 3, 4 and 5. 
For Case 4
we use another interval to 
describe equalities on $w(Q^n)$ in Section 6,
but it is closely related to the intervals and the rectangle introduced here. 
Although we do not use blow-ups in the proofs of our main theorems,
they are useful to explain our results in both the previous and this papers. 
Let 
\[
\alpha = \dfrac{\gamma}{\delta - d}
\]
if $\delta \neq d$, and
we assume that $a_{\delta} = 1$ and $b_{\gamma d} = 1$ for simplicity. 

\subsection{Interval of weights and Blow-ups for Case 2}

Let $s > 1$,
\[
\delta \leq T_{s-1}, \
(\gamma, d) = (n_s, m_s)
\text{ and } 
l_1 = \frac{n_s - n_{s-1}}{m_{s-1} - m_s}.
\]
Note that $\gamma > 0$ and
$j \geq d$ for any $(i,j)$ such that $b_{ij} \neq 0$ by the setting.

We define the interval $\mathcal{I}_f$ as 
\[
\mathcal{I}_f = 
\left\{ \ l > 0 \ | 
\begin{array}{lcr}
l \delta \leq \gamma + ld \leq i + l j
\text{ for any $i$ and $j$ such that } b_{ij} \neq 0
\end{array}
\right\}.
\]
If $\delta > d$, then 
\[
\mathcal{I}_f 
=
\left[
\max_{i,j} \left\{ \dfrac{\gamma - i}{j - d} \ \Big| \ b_{ij} \neq 0 \text{ and } j > d \right\},
\dfrac{\gamma}{\delta - d}
\right].
\]
Since the ratio $(d - j)/(\gamma - i)$ is the slope of 
the line passing through $(\gamma, d)$ and $(i,j)$,
we can take the maximum over $(i,j)$ at $(n_{s-1}, m_{s-1})$:
\[
\mathcal{I}_f 
=
\left[
\max_{1 \leq j < s} \left\{ \dfrac{\gamma - n_j}{m_j - d} \right\},
\dfrac{\gamma}{\delta - d}
\right]
=
\left[
\frac{\gamma - n_{s-1}}{m_{s-1} - d},
\dfrac{\gamma}{\delta - d}
\right]
=
\left[
l_1,
\alpha
\right],
\]
which is mapped to $[\delta, T_{s-1}]$
by the transformation $l \to l^{-1} \gamma + d$.
Therefore,
$\mathcal{I}_f $ can be identified with the set of the lines
passing through $(\gamma, d)$ 
whose slopes are in $[- \alpha^{-1}, -l_1^{-1}]$
or, equivalently,
whose $y$-intercepts are in $[\delta, T_{s-1}]$.
If $\delta \leq d$, then
the inequality $l \delta \leq \gamma + ld$ is trivial and so
$\mathcal{I}_f =[ l_1, \infty )$.
In particular, 
$\min \mathcal{I}_f = l_1$.


Assuming that $l$ in $\mathcal{I}_f$ is an integer,
we explain 
our previous results in terms of a blow-up. 
Let $\pi_1 (z,c) = (z, z^{l} c)$
and $\tilde{f} = \pi_1^{-1} \circ f \circ \pi_1$.
Note that $\pi_1$ is the $l$-th composition of the blow-up $(z,c) \to (z,zc)$. 
Then we have
\begin{align*}
\tilde{f} (z,c) 
&= (p(z), \tilde{q}(z,c))
= \left( p(z), \ \dfrac{q(z,z^l c)}{p(z)^l} \right) \\
&= \left( z^{\delta} (1 + o(z)), \ 
\sum b_{ij} z^{i + lj - l \delta} c^{j} (1 + o(z)) \right).
\end{align*}
Let $\tilde{\gamma} = \gamma + l d - l \delta$,
$\tilde{i} = i + l j - l \delta$
and $\tilde{n}_j = n_j + l m_j - l \delta$.
Then 
the Newton polygon $N(\tilde{q})$ of $\tilde{q}$ has just one vertex $(\tilde{\gamma}, d)$:
$N(\tilde{q}) = D(\tilde{\gamma}, d)$.

\begin{lemma} \label{blow-up lem for case 2} 
It follows that
$0 \leq \tilde{\gamma} \leq \tilde{i}$ 
for any $l$ in $\mathcal{I}_f$ and for any $(i,j)$ such that $b_{ij} \neq 0$.
In particular,
$0 \leq \tilde{\gamma} \leq \tilde{n}_j$ 
for any $l$ in $\mathcal{I}_f$ and for any $1 \leq j \leq s$.
More precisely,
$\tilde{\gamma} = \tilde{n}_{s-1}$ and
$\tilde{\gamma} < \tilde{n}_{j}$ for any $j \neq s-1$, $s$ if $l = l_1$, and
$\tilde{\gamma} < \tilde{n}_j$ for any $j \neq s$ if $l > l_1$.
Moreover,
$\tilde{\gamma} > 0$ if $l < \alpha$, and 
$\tilde{\gamma} = 0$ if $l = \alpha$.  
\end{lemma}

\begin{remark} \label{A_1}
The blow-up of $f$ can be transfered to the affine transformation of $N(q)$.
The affine transformation
\[
A_1
\begin{pmatrix} i \\ j \end{pmatrix} =
\begin{pmatrix} i + l_1 j - l_1 \delta \\ j \end{pmatrix} =
\begin{pmatrix} 1 & l_1 \\ 0 & 1 \end{pmatrix}
\begin{pmatrix} i \\ j \end{pmatrix} - \begin{pmatrix} l_1 \delta \\ 0 \end{pmatrix}
\]
maps the basis $\{ (1,0), (-l_1,1) \}$ to $\{ (1,0), (0,1) \}$.
In other words, 
$A_1$ maps a horizontal line and the line $L_{s-1}$ with slope $-l_1^{-1}$
to the same horizontal line and a vertical line.
\end{remark}

\begin{proposition} \label{} 
If $l$ in $\mathcal{I}_f$ is an integer, 
then $\tilde{f}$ is well-defined, holomorphic and skew product  
on a neighborhood of the origin. 
More precisely,
\[
\tilde{f} (z,c) = \big( z^{\delta} (1 + o(z)), \ 
z^{\tilde{\gamma}} c^d (1 + o(z, c)) \big), 
\]
and it has a fixed point at the origin if $d > 0$.
Moreover,
if $\delta \geq 2$ and $d \geq 2$
or if $\delta \geq 2$, $d=1$ and $\delta < T_{s-1}$,
then
the fixed point is superattracting.
\end{proposition}

Therefore,
if $l$ is an integer and $\tilde{f}$ is superattracting,
then it is easy to construct the B\"{o}ttcher coordinate for $\tilde{f}$
on a neighborhood of the origin
that conjugates $\tilde{f}$ to
$(z,c) \to (z^{\delta}, z^{\tilde{\gamma}} c^d)$,
because $\tilde{f}$ is a holomorphic skew product in Case 1.
Consequently,
we obtain the B\"{o}ttcher coordinate for $f$ on $U^l$
that conjugates $f$ to $f_0$,
where
$f_0(z,w) = (z^{\delta}, z^{\gamma} w^d)$
and $U^l = \{ |z| < r, |w| < r |z|^{l} \}$.
Actually,
we can construct the B\"{o}ttcher coordinate for $f$ on $U^l$ directly
even if $l$ in $\mathcal{I}_f$ is not an integer. 

\begin{theorem}[\cite{u}] \label{Bottcher for case 2} 
If $\delta \geq 2$ and $d \geq 2$
or if $\delta \geq 2$, $d=1$ and $\delta < T_{s-1}$, 
then for any $l$ in $\mathcal{I}_f$,
$f$ preserves $U^l$ and 
there is a biholomorphic map defined on $U^l$
that conjugates $f$ to $f_0$
for small $r$.
\end{theorem}

Note that $U = U^{l_1}$, 
which is the largest region among $U^l$ for any $l$ in $\mathcal{I}_f$.

\begin{remark} \label{}  
Even if $l$ is rational,
we can lift $f$ to a holomorphic skew product.
In fact,
let $\pi_1 (\mathsf{z}, c) = (\mathsf{z}^r, \mathsf{z}^s c)$
and $\tilde{f} = \pi_1^{-1} \circ f \circ \pi_1$,
where $s/r = l$.
Then 
\[
\tilde{f} (\mathsf{z},c) 
= \left( \mathsf{z}^{\delta} (1 + o(\mathsf{z})), \ 
\sum b_{ij} \mathsf{z}^{ri + sj - s \delta} c^{j} (1 + o(\mathsf{z})) \right)
\]
and it is well-defined, holomorphic and skew product. 
\end{remark}

\subsection{Interval of weights and Blow-ups for Case 3}

Let $s > 1$,
\[
T_1 \leq \delta, \
(\gamma, d) = (n_1, m_1) \text{ and }
l_2 = \frac{n_2 - n_1}{m_1 - m_2}. 
\]
Note that $\delta \geq d > 0$ and 
$i \geq \gamma$ for any $(i,j)$ such that $b_{ij} \neq 0$ by the setting.

We define the interval $\mathcal{I}_f$ as
\[
\mathcal{I}_f = 
\left\{ \ l > 0 \ \Big| 
\begin{array}{lcr}
\gamma + ld \leq i + lj
\text{ and }
\gamma + ld \leq l \delta \\ 
\text{for any $i$ and $j$ such that } b_{ij} \neq 0
\end{array}
\right\}. 
\]
If $\gamma > 0$, then $\delta > d$ and 
\begin{align*}
\mathcal{I}_f 
&= 
\left[
\dfrac{\gamma}{\delta - d},
\min_{i,j} \left\{ \dfrac{i - \gamma}{d - j} \ \Big| 
\begin{array}{lcr}
b_{ij} \neq 0 \text{ and } j < d 
\end{array}
\right\}
\right] \\
&= 
\left[
\dfrac{\gamma}{\delta - d},
\min_{1 < j \leq s} \left\{ \dfrac{n_j - \gamma}{d - m_j} \right\}
\right]
=
\left[
\dfrac{\gamma}{\delta - d},
\frac{n_2 - \gamma}{d - m_2}
\right]
=
\left[
\alpha,
l_2
\right],
\end{align*}
which is mapped to $[T_1, \delta]$
by the transformation $l \to l^{-1} \gamma + d$.
Therefore,
$\mathcal{I}_f $ can be identified with the set of the lines
passing through $(\gamma, d)$ 
whose slopes are in $[-l_2^{-1}, - \alpha^{-1}]$
or, equivalently,
whose $y$-intercepts are in $[\delta, T_{1}]$.
If $\gamma = 0$, then
the inequality $\gamma + ld \leq l \delta$ is trivial since $d \leq \delta$
and so $\mathcal{I}_f = (0, l_2]$.
In particular,
$\max \mathcal{I}_f = l_2$.


Assuming that $l^{-1}$ is an integer for $l$ in $\mathcal{I}_f$,
we explain 
our previous results in terms of a blow-up. 
Let $\pi_2 (t,w) = (t w^{l^{-1}}, w)$ and $\tilde{f} = \pi_2^{-1} \circ f \circ \pi_2$.
Note that $\pi_2$ is the $l^{-1}$-th composition of the blow-up $(t,w) \to (tw,w)$.
Then we have
\[
\tilde{f} (t,w) 
= (\tilde{p} (t,w), \tilde{q} (t,w))
= \left( \dfrac{p(tw^{l^{-1}})}{q(tw^{l^{-1}}, w)^{l^{-1}}}, \ q(tw^{l^{-1}}, w) \right).
\]
Let $\tilde{d} = l^{-1} \gamma + d$ and $\tilde{j} = l^{-1} i + j$.
Then $\tilde{q}(t,w) = \sum b_{ij} t^{i} w^{l^{-1} i + j} = \sum b_{ij} t^{i} w^{\tilde{j}}$
and the Newton polygon $N(\tilde{q})$ of $\tilde{q}$ 
has just one vertex $(\gamma, \tilde{d})$:
$N(\tilde{q}) = D(\gamma, \tilde{d})$.

\begin{lemma} \label{blow-up lem for case 3} 
It follows that
$d < \tilde{d} \leq \tilde{j}$ 
for any $l$ in $\mathcal{I}_f$ and for any $(i,j)$ such that $b_{ij} \neq 0$.
\end{lemma}

\begin{remark}
The linear transformation 
\[
A_2 
\begin{pmatrix} i \\ j \end{pmatrix} =
\begin{pmatrix} i \\ l_2^{-1} i + j \end{pmatrix} =
\begin{pmatrix} 1 & 0 \\ l_2^{-1} & 1 \end{pmatrix}
\begin{pmatrix} i \\ j \end{pmatrix}
\]
maps the basis $\{ (1,-l_2^{-1}), (0,1) \}$ to $\{ (1,0), (0,1) \}$.
In other words, 
$A_2$ maps the line $L_{1}$ with slope $-l_2^{-1}$ and a vertical line 
to a horizontal line and the same vertical line.
\end{remark}

Hence $\tilde{q} (t,w) = t^{\gamma} w^{l^{-1} \gamma + d} (1 + o(t,w)) 
                            = t^{\gamma} w^{\tilde{d}} (1 + o(t,w))$ and so
\[
\tilde{p} (t,w) 
= \dfrac{ \{ tw^{l^{-1}} \}^{\delta} (1 + o(t,w)) }{ \{ t^{\gamma} w^{\tilde{d}} (1 + o(t,w)) \}^{l^{-1}} }
= t^{\delta - l^{-1} \gamma} w^{l^{-1} (\delta - \tilde{d})} (1 + o(t,w)).
\]
Note that $\tilde{d} \leq \delta$ since $\gamma + ld \leq l \delta$.
Hence $\delta - l^{-1} \gamma \geq d > 0$ and $\delta - \tilde{d} \geq 0$.

\begin{proposition} \label{blow-up prop for case 3} 
If $l^{-1}$ is an integer for $l$ in $\mathcal{I}_f$,
then $\tilde{f}$ is well-defined and holomorphic 
on a neighborhood of the origin. 
More precisely,
\[
\tilde{f} (t,w) = 
\big( t^{\delta - l^{-1} \gamma} w^{l^{-1} (\delta - \tilde{d})} (1 + o(t,w)),
\ t^{\gamma} w^{\tilde{d}} (1 + o(t,w)) \big),
\]
and it has a fixed point at the origin.
Moreover,
if $\delta \geq 2$ and $d \geq 2$
or if $\delta \geq 2$, $d = 1$ and $\delta > T_1$,
then
the fixed point is superattracting.
\end{proposition}

Although $\tilde{f}$ is not skew product,
it is a perturbation of a monomial map near the origin.
Therefore,
if $l^{-1}$ is an integer and $\tilde{f}$ is superattracting,
then we can construct the B\"{o}ttcher coordinate for $\tilde{f}$,
which induces the B\"{o}ttcher coordinate for $f$ on $U^l$,
where
$U^l = \{ |z| < r|w|^l, |w| < r \}$.
Actually,
we can construct the B\"{o}ttcher coordinate for $f$ on $U^l$ directly
for any $l$ in $\mathcal{I}_f$.

\begin{theorem}[\cite{u}] \label{Bottcher for case 3} 
If $\delta \geq 2$ and $d \geq 2$
or if $\delta \geq 2$, $d=1$ and $\delta > T_{1}$, 
then for any $l$ in $\mathcal{I}_f$,
$f$ preserves $U^l$ and
there is a biholomorphic map defined on $U^l$
that conjugates $f$ to $f_0$
for small $r$.
\end{theorem}

Note that $U = U^{l_2}$,
which is the largest region among $U^l$ for any $l$ in $\mathcal{I}_f$.

\begin{remark} \label{} 
If $l^{-1}$ is rational, 
then a similar lift is well-defined
only for limited cases.
Let $\pi_2 (t, \mathsf{w}) = (t \mathsf{w}^r, \mathsf{w}^s)$ 
and $\tilde{f} = \pi_2^{-1} \circ f \circ \pi_2$,
where $s/r = l$. 
Then $\tilde{f}$ is well-defined if $\gamma/s$ is an integer.
\end{remark}

\subsection{Intervals of weights and Blow-ups for Case 4}

Let $s > 2$,
\[
T_k \leq \delta \leq T_{k-1} 
\text{ for some }
2 \leq k \leq s-1,
\] 
\[
(\gamma, d) = (n_k, m_k), \
l_1 = \frac{n_k - n_{k-1}}{m_{k-1} - m_k} \text{ and } 
l_1 + l_2= \frac{n_{k+1} - n_k}{m_k - m_{k+1}}.
\] 
Note that $\gamma d > 0$ and $\delta > d$ by the setting.

We define the interval $\mathcal{I}_f^1$ as
\[
\mathcal{I}_f^1 = 
\left\{ \ l_{(1)} > 0 \ \Bigg| 
\begin{array}{lcr}
\gamma + l_{(1)} d \leq n_{j} + l_{(1)} m_{j} \text{ for any } j \leq k-1 \\ 
\gamma + l_{(1)} d  <   n_{j} + l_{(1)} m_{j} \text{ for any } j \geq k+1 \\ 
l_{(1)} \delta \leq \gamma + l_{(1)} d 
\end{array}
\right\},
\]
the interval $\mathcal{I}_f^2$ associated with $l_{(1)}$ in $\mathcal{I}_f^1$ as
\[
\mathcal{I}_f^2  = \mathcal{I}_f^2 (l_{(1)}) = 
\left\{ \ l_{(2)} > 0 \ \Big| 
\begin{array}{lcr}
\tilde{\gamma} + l_{(2)} d \leq \tilde{i} + l_{(2)} j
\text{ and }
\tilde{\gamma} + l_{(2)} d \leq l_{(2)} \delta  \\
\text{for any $i$ and $j$ such that } b_{ij} \neq 0 
\end{array}
\right\},
\]
where 
$\tilde{\gamma} = \gamma + l_{(1)} d - l_{(1)} \delta$
and $\tilde{i} = i + l_{(1)} j - l_{(1)} \delta$,
and the rectangle $\mathcal{I}_f$ as
\[
\mathcal{I}_f =
\left\{ 
(l_{(1)}, l_{(1)} + l_{(2)}) \ | \ l_{(1)} \in \mathcal{I}_f^1, l_{(2)} \in \mathcal{I}_f^2 (l_{(1)}) 
\right\}.
\]

Let us calculate the intervals and rectangle more practically.
Note that $\alpha > 0$ since $\delta > d$ and $\gamma > 0$ by the setting.
Since $n_j < \gamma$ and $m_j > d$ for any $j \leq k-1$,
and $n_j > \gamma$ and $m_j < d$ for any $j \geq k+1$,  
\begin{align*}
\mathcal{I}_f^1 &= 
\left[
\max_{j \leq k-1}
\left\{ 
\dfrac{\gamma - n_j}{m_j - d} 
\right\},
\min_{j \geq k+1}
\left\{ 
\dfrac{n_j - \gamma}{d - m_j} 
\right\}
\right)
\cap
\left(
0,
\dfrac{\gamma}{\delta - d}
\right] \\
&= 
\left[
\dfrac{\gamma - n_{k-1}}{m_{k-1} - d},
\dfrac{n_{k+1} - \gamma}{d - m_{k+1}} 
\right)
\cap
\left(
0,
\dfrac{\gamma}{\delta - d}
\right]
= [ l_1, l_1 + l_2 ) \cap ( 0, \alpha ].
\end{align*}
In particular,
$\min \mathcal{I}_f^1 = l_1$ and, as a remark,
\[
\mathcal{I}_f^1 = 
\left\{ \ l_{(1)} > 0 \ \Bigg| 
\begin{array}{lcr}
\gamma + l_{(1)} d \leq n_{k-1} + l_{(1)} m_{k-1} \\
\gamma + l_{(1)} d < n_{k+1} + l_{(1)} m_{k+1} \\
l_{(1)} \delta \leq \gamma + l_{(1)} d 
\end{array}
\right\}.
\]
On the other hand,
\begin{align*}
\mathcal{I}_f^2 
&=
\left[
\dfrac{\tilde{\gamma}}{\delta - d},
\dfrac{\tilde{n}_{k+1} - \tilde{\gamma}}{d - m_{k+1}}  
\right]
\cap \mathbb{R}_{>0}
=
\left[
\dfrac{\gamma}{\delta - d} - l_{(1)},
\dfrac{n_{k+1} - \gamma}{d - m_{k+1}} - l_{(1)}  
\right]
\cap \mathbb{R}_{>0} \\
&= [ \alpha - l_{(1)}, l_1 + l_2 - l_{(1)} ]
\cap \mathbb{R}_{>0}.
\end{align*}
If $T_k < \delta = T_{k-1}$, 
then it follows from the inequality $l_1 = \alpha < l_1 + l_2$ that
\[
\mathcal{I}_f^1 = \{ l_1 \}, \ 
\mathcal{I}_f^2 = (0, l_2]  
\text{ and so }
\mathcal{I}_f 
= \{ l_1 \} \times ( l_1, l_1 + l_2 ]. 
\]
If $T_k < \delta < T_{k-1}$, 
then it follows from the inequality $l_1 < \alpha < l_1 + l_2$ that
\[
\mathcal{I}_f^1 = [ l_1, \alpha ], \ 
\mathcal{I}_f^2 = 
  \begin{cases} 
    [ \alpha - l_{(1)}, l_1 + l_2 - l_{(1)} ] & \text{ if } l_{(1)} < \alpha, \\
    (0, l_1 + l_2 - \alpha ] & \text{ if } l_{(1)} = \alpha
  \end{cases}
\]
\[
\text{and so }
\mathcal{I}_f 
= [ l_1, \alpha ] \times [ \alpha, l_1 + l_2 ] - \{ (\alpha, \alpha) \}.
\]
If $T_{k} = \delta < T_{k-1}$, then it follows from the inequality $l_1 < \alpha = l_1 + l_2$ that
\[
\mathcal{I}_f^1 = [ l_1, l_1 + l_2 ), \ 
\mathcal{I}_f^2 = \{ l_1 + l_2 - l_{(1)} \}
\text{ and so }
\mathcal{I}_f 
= [ l_1, l_1 + l_2 ) \times \{ l_1 + l_2 \}.
\]
In particular, 
$\min \mathcal{I}_f^1 = l_1$ and
$\max \{ l_{(1)} + l_{(2)} \ | \ l_{(1)} \in \mathcal{I}_f^1, l_{(2)} \in \mathcal{I}_f^2 (l_{(1)}) \} = l_1 + l_2$.


Assuming that $l_1$ and $l_2^{-1}$ are integers,
we explain
our previous results in terms of blow-ups. 
Although the same arguments hold for 
any $l_{(1)}$ in $\mathcal{I}_f^1$ and $l_{(2)}$ in $\mathcal{I}_f^2$,
the case $l_{(1)} = l_1$ and $l_{(2)} = l_2$ is most important 
for the study in the previous and this papers, and so
we fix $l_1$ and $l_2$ for simplicity.
The strategy is to combine the blow-ups in Cases 2 and 3.
We first blow-up $f$ to $\tilde{f}_1$
by $\pi_1$ as in Case 2.
It then turns out that $\tilde{f}_1$ is a holomorphic skew product in Case 3.
We next blow-up $\tilde{f}_1$ to $\tilde{f}_2$
by $\pi_2$ as in Case 3.
The map $\tilde{f}_2$ is a perturbation of a monomial map near the origin,
and we obtain the B\"{o}ttcher coordinate for $f$ on $U$ 
if $\tilde{f}_2$ is superattracting. 


Let us explain the first blow-up precisely.
Let $\pi_1 (z,c) = (z, z^{l_1} c)$ and 
$\tilde{f}_1 = \pi_1^{-1} \circ f \circ \pi_1$ as in Case 2.
Then 
\begin{align*}
\tilde{f}_1 (z,c) 
&= (\tilde{p}_1 (z), \tilde{q}_1 (z,c))
= \left( p(z), \ \dfrac{q(z,z^{l_1} c)}{p(z)^{l_1}} \right) \\
&= \left( z^{\delta} (1 + o(z)), \ 
\sum b_{ij} z^{ i + l_1 j - l_1 \delta} c^{j} (1 + o(z, c)) \right).
\end{align*}

\begin{proposition} \label{} 
If $l_1$ is an integer, 
then $\tilde{f}_1$ is well-defined, holomorphic and skew product 
on a neighborhood of the origin,
which has a fixed point at the origin.  
Moreover,
if $\delta \geq 2$ and $d \geq 2$
or if $\delta \geq 2$, $d = 1$ and $T_k < \delta < T_{k-1}$,
then
the fixed point is superattracting.
\end{proposition}

Let $\tilde{\gamma} = \gamma + l_1 d - l_1 \delta$,
$\tilde{i} = i + l_1 j - l_1 \delta$
and $\tilde{n}_j = n_j + l_1 m_j - l_1 \delta$
as in Case 2.
Then $0 \leq \tilde{\gamma} \leq \tilde{n}_j$ for any $1 \leq j \leq s$.
In particular,
$(\tilde{\gamma}, d)$ is minimum in the sense that
$\tilde{\gamma} \leq \tilde{i}$, and $d \leq j$ if $\tilde{\gamma} = \tilde{i}$.
Hence
$(\tilde{\gamma}, d)$ is the vertex of the Newton polygon $N(\tilde{q}_1)$
whose $x$-coordinate is minimum.
However,
$N(\tilde{q}_1)$ has other vertices such as $(\tilde{n}_{k+1}, m_{k+1})$.
Hence the situation resembles that of Case 3.

We illustrate that $\tilde{f}_1$ is actually in Case 3. 
Recall that $L_{k}$ is the line passing through 
the vertices $(\gamma, d)$ and $(n_{k+1}, m_{k+1})$,
and $T_{k}$ is the $y$-intercept of $L_{k}$.
The slope of $L_{k}$ is $-(l_1 + l_2)^{-1}$ 
and so $T_{k} = (l_1 + l_2)^{-1} \gamma + d$.
Let $\tilde{L}_{k}$ be the line passing through  
the vertices $(\tilde{\gamma}, d)$ and $(\tilde{n}_{k+1}, m_{k+1})$,
and let $\tilde{T}_{k}$ be the $y$-intercept of $\tilde{L}_{k}$.
Then the slope of $\tilde{L}_{k}$ is $-l_2^{-1}$
and so $\tilde{T}_{k} = l_2^{-1} \tilde{\gamma} + d$
because the affine transformation $A_1$ in Remark \ref{A_1}  
maps 
the basis $\{ (1,0), (-(l_1 + l_2),1) \}$ to $\{ (1,0), (-l_2,1) \}$.
Moreover,
$\tilde{T}_k \leq \delta$ since $T_k \leq \delta$. 
More precisely,
$\tilde{T}_k < \delta$ if $T_k < \delta$,
and $\tilde{T}_k = \delta$ if $T_k = \delta$.

\begin{proposition} \label{} 
If $l_1$ is an integer, 
then $\tilde{f}_1$ is a holomorphic skew product in Case 3. 
\end{proposition}


Let us next explain the second blow-up precisely.
Let $\pi_2 (t,c) = (t c^{l_2^{-1}}, c)$ and $\tilde{f}_2 = \pi_2^{-1} \circ \tilde{f}_1 \circ \pi_2$ as in Case 3.
Then 
\[
\tilde{f}_2 (t,c) 
= (\tilde{p}_2 (t,c), \tilde{q}_2 (t,c))
= \left( \dfrac{\tilde{p}_1 (tc^{l_2^{-1}})}{\tilde{q}_1 (tc^{l_2^{-1}},c)^{l_2^{-1}}}, \ \tilde{q}_1 (tc^{l_2^{-1}},c) \right).
\]
Let $\tilde{d} = l_2^{-1} \tilde{\gamma} + d$
and $\tilde{j} = l_2^{-1} \tilde{i} + j$
as in Case 3.
Then
$\tilde{d} \leq \tilde{j}$ for any $(i,j)$ such that $b_{ij} \neq 0$, and $\tilde{d} \leq \delta$.
In particular,
the minimality of $(\tilde{\gamma}, \tilde{d})$ follows.

\begin{lemma} \label{blow-up lem for case 4}
It follows that
$0 \leq \tilde{\gamma} \leq \tilde{i}$ and 
$d \leq \tilde{d} \leq \tilde{j}$ for any $(i,j)$ such that $b_{ij} \neq 0$.
\end{lemma}

Hence
the Newton polygon $N(\tilde{q}_2)$ of $\tilde{q}_2$ has just one vertex $(\tilde{\gamma}, \tilde{d})$:
$N(\tilde{q}_2) = D(\tilde{\gamma}, \tilde{d})$.

\begin{proposition} \label{blow-up prop for case 4} 
If $l_1$ and $l_2^{-1}$ are integers, 
then $\tilde{f}_2$ is well-defined and holomorphic 
on a neighborhood of the origin.  
More precisely,
\[
\tilde{f}_2 (t,c) =
\big( t^{\delta - l_2^{-1} \tilde{\gamma}} c^{l_2^{-1} (\delta - \tilde{d})} ( 1 +  o(t,c) ), 
\ t^{\tilde{\gamma}} c^{\tilde{d}} ( 1 + o(t,c) ) \big),
\]
and it has a fixed point at the origin.
Moreover,
if $\delta \geq 2$ and $d \geq 2$
or if $\delta \geq 2$, $d = 1$ and $T_k < \delta < T_{k-1}$,
then the fixed point is superattracting.
\end{proposition}

Therefore,
if $l_1$ and $l_2^{-1}$ are integers and 
if $\tilde{f}_2$ is superattracting,
then we can construct the B\"{o}ttcher coordinate for $\tilde{f}_2$
on a neighborhood of the origin, 
which induces that for $\tilde{f}_1$ 
on an open wedge and  
that for $f$ 
on $U$.
Actually,
we can construct the B\"{o}ttcher coordinate for $f$ on $U$ directly
even if $l_1$ and $l_2^{-1}$ are not integers and,
moreover, 
the B\"{o}ttcher coordinate for $f$ on $U^{l_{(1)}, l_{(2)}}$ directly
for any $l_{(1)}$ in $\mathcal{I}_f^1$ and $l_{(2)}$ in $\mathcal{I}_f^2$,
where
$U^{l_{(1)}, l_{(2)}} = \{ r^{-l_{(2)}} |z|^{l_{(1)} + l_{(2)}} < |w| < r|z|^{l_{(1)}} \}$.

\begin{theorem}[\cite{u}] \label{Bottcher for case 4} 
If $\delta \geq 2$ and $d \geq 2$
or if $\delta \geq 2$, $d=1$ and $T_{k} < \delta < T_{k-1}$, 
then for any $l_{(1)}$ in $\mathcal{I}_f^1$
and $l_{(2)}$ in $\mathcal{I}_f^2$,
$f$ preserves $U^{l_{(1)}, l_{(2)}}$ and 
there is a biholomorphic map defined on $U^{l_{(1)}, l_{(2)}}$
that conjugates $f$ to $f_0$
for small $r$.
\end{theorem}

Note that $U = U^{l_1, l_2}$,
which is the largest region among $U^{l_{(1)}, l_{(2)}}$
for any $l_{(1)}$ in $\mathcal{I}_f^1$ and $l_{(2)}$ in $\mathcal{I}_f^2$.


\begin{remark}
The affine transformation
\[
A 
\begin{pmatrix} i \\ j \end{pmatrix} =
\begin{pmatrix} 1 & 0 \\ l_2^{-1} & 1 \end{pmatrix}
\left\{  
\begin{pmatrix} 1 & l_1 \\ 0 & 1 \end{pmatrix}
\begin{pmatrix} i \\ j \end{pmatrix} - 
\begin{pmatrix} l_1 \delta \\ 0 \end{pmatrix}
\right\}
\]
is the composition of two affine transformations
\[
A_1
\begin{pmatrix} i \\ j \end{pmatrix} =
\begin{pmatrix} i + l_1 j - l_1 \delta \\ j \end{pmatrix} 
\text{ and }
A_2
\begin{pmatrix} i \\ j \end{pmatrix} =
\begin{pmatrix} i \\ l_2^{-1} i + j \end{pmatrix}. 
\]
The transformation $A_1$ maps 
the basis $\{ (1,-(l_1 +l_2)^{-1}), (-l_1,1) \}$ to $\{ (1,-l_2^{-1}), (0,1) \}$. 
In other words, 
it maps the line $L_{k}$ with slope $-(l_1 +l_2)^{-1}$ and the line $L_{k-1}$ with slope $-l_1^{-1}$,
which intersect at $(\gamma, d)$,
to the line $\tilde{L}_{k}$ with slope $-l_2^{-1}$ and the vertical line,
which intersect at $(\tilde{\gamma}, d)$.
The transformation $A_2$ maps 
the basis $\{ (1,-l_2^{-1}),  (0,1) \}$ to $ \{ (1,0), (0,1) \}$.
In other words, 
it maps the line $\tilde{L}_{k}$ and the vertical line,
which intersect at $(\tilde{\gamma}, d)$, 
to the horizontal line and the vertical line,
which intersect at $(\tilde{\gamma}, \tilde{d})$.
Consequently, 
$A$ maps the lines $L_{k}$ and $L_{k-1}$
to the horizontal and vertical lines.
\end{remark}

\section{Shape of Newton polygon of $Q^n$ for Case 2}

In this section we deal with Case 2
under the condition $d > 0$.
Let $s > 1$,
\[
\delta \leq T_{s-1}, \
(\gamma, d) = (n_s, m_s)
\text{ and } 
l_1 = \frac{n_s - n_{s-1}}{m_{s-1} - m_s}.
\]
Note that $\gamma > 0$ by the setting.
We first give a summary of the results in Section 3.1
and illustrate some of them in terms of the blow-ups in Section 3.2. 
Preparing two lemmas in Section 3.3,
we prove Theorem \ref{main thm of weights when d>0} below in Section 3.4.
Finally,
we show the existence of the vertex of $N(Q^n)$
that is previous to $(\gamma_n, d^n)$
and specify it in Section 3.5,
which induces Theorem \ref{main thm of attr rates when d>0} below.

\subsection{Summary of results}

We first show the following equalities.

\begin{theorem} \label{main thm of weights when d>0}
Let $d > 0$. Then 
$Q^n$ contains the term $z^{\gamma_n} w^{d^n}$ and
$w(Q^n) = \gamma_n + l d^n$
for any $n \geq 1$ and
for any $l$ in $\mathcal{I}_f$,
where $\mathcal{I}_f = [l_1, \alpha]$ or $\mathcal{I}_f = [l_1, \infty)$
if $\delta > d$ or $\delta \leq d$.
\end{theorem}

Moreover, 
$(\gamma_n, d^n)$ is the vertex of $N(Q^n)$
whose $y$-coordinate is minimum.
Hence the theorem above implies the following corollary. 

\begin{corollary} \label{main cor on attr rates when d>0}
Let $d > 0$.
Then it follows for any $n \geq 1$ that
\begin{enumerate}
\item[$(1)$] $c(Q^n) = \gamma_n + d^n$ 
if $l_1 \leq 1$, and
\item[$(2)$] $l_1^{-1} \gamma_n + d^n \leq c(Q^n) \leq \gamma_n + d^n$
if $l_1 > 1$.
\end{enumerate}
\end{corollary}

Furthermore,
we can specify the vertex of $N(Q^n)$ that is previous to $(\gamma_n, d^n)$,
which exists,
and improve the corollary above as follows.

\begin{theorem} \label{main thm of attr rates when d>0}
Let $d > 0$.
Then it follows for any $n \geq 1$ that
\begin{enumerate}
\item[$(1)$] $c(Q^n) = \gamma_n + d^n$
if $l_1 \leq 1$, and
\item[$(2)$] $l_1^{-1} \gamma_n + d^n \leq c(Q^n) < \gamma_n + d^n$
if $l_1 > 1$.
\end{enumerate}
More precisely,
it follows for any $n \geq 1$ that
\begin{enumerate}
\item[$(3)$] $l_1^{-1} \gamma_n + d^n < c(Q^n)$
if $l_1 > 1$ and if $n_1 > 0$ or $s > 2$.
\end{enumerate}
Let $d > 0$, $l_1 > 1$, $n_1 = 0$ and $s = 2$.
Then
\begin{enumerate}
\item[$(4)$] $l_1^{-1} \gamma + d = c(q)$ and $l_1^{-1} \gamma_n + d^n < c(Q^n)$ 
for any $n \geq 2$ if $\delta < T_{s-1}$, and
\item[$(5)$] $l_1^{-1} \gamma_n + d^n = c(Q^n)$
for any $n \geq 1$ if $\delta = T_{s-1}$.
\end{enumerate}
\end{theorem}

\subsection{Illustration of results in terms of blow-ups}

Assuming that $l_1$ is an integer,
we illustrate Theorem \ref{main thm of weights when d>0} and 
Corollary \ref{main cor on attr rates when d>0} in terms of the blow-up. 
Let $\pi_1 (z,c) = (z, z^{l_1} c)$
and $\tilde{f} = \pi_1^{-1} \circ f \circ \pi_1$
as in Section 2.1.
By Lemma \ref{blow-up lem for case 2},
$N(\tilde{q}) = D(\tilde{\gamma}, d)$,
where $\tilde{\gamma} = \gamma + l_1 d - l_1 \delta$.
Hence $\tilde{f}$ is a skew product in Case 1.
If $d>0$, 
then $\tilde{Q}^n$ contains the term $z^{\tilde{\gamma}_n} w^{d^n}$ and
$N(\tilde{Q}^n) = D(\tilde{\gamma}_n, d^n)$,
where $\tilde{\gamma}_n = \gamma_n + l_1 d^n - l_1 \delta^n$.
This implies that 
$Q^n$ contains the term $z^{\gamma_n} w^{d^n}$, and that
$N(Q^n)$ is included in the upper-right region
that is surrounded by the line with slope $-l_1^{-1}$
and the horizontal line,
which intersect at $(\gamma_n, d^n)$. 
Therefore,
we obtain Theorem \ref{main thm of weights when d>0}
and Corollary \ref{main cor on attr rates when d>0} when $l_1$ is an integer,
because
$w(Q^n)$ and $c(Q^n)$ are the minimum $x$-intercepts of 
the lines with slopes $-l^{-1}$ and $-1$ 
that intersect $N(Q^n)$.

\subsection{Preliminary lemmas: Affine dynamics on the interval}

The following affine function $R$ on the real line plays an important role:
\[
R(l) = \dfrac{\gamma + l d}{\delta}.
\]
If $\delta \neq d > 0$, then 
$R$ is a contracting or expanding function around the fixed point $\alpha$:
\[
R(l) = \dfrac{d}{\delta} (l - \alpha) + \alpha
\text{ and so }
R^n (l) = \left( \dfrac{d}{\delta} \right)^n (l - \alpha) + \alpha.
\]
More precisely,
$\alpha$ is attracting and $\alpha > 0$ if $\delta > d$, and
$\alpha$ is repelling and $\alpha < 0$ if $\delta < d$.
On the other hand,
if $\delta = d$, then $R$ is a translation:
\[
R(l) = l + \dfrac{\gamma}{\delta}
\text{ and so }
R^n (l) = l + \dfrac{\gamma}{\delta} n.
\]
Since $\mathcal{I}_f = [l_1, \alpha]$ or $\mathcal{I}_f = [l_1, \infty)$
if $\delta > d$ or $\delta \leq d$,
we can conclude as follows.

\begin{lemma} \label{dyn of R when d>0}
For any $l$ in $\mathcal{I}_f$,
it follows that $R^n (l)$ belongs to $\mathcal{I}_f$ for any $n \geq 1$ and 
the sequence $\{ R^n (l) \}_{n \geq 1}$ is increasing. 
More precisely, 
$R^n (l) \to \alpha$ or $R^n (l) \to \infty$ as $n \to \infty$
if $\delta > d > 0$ or $\delta \leq d$.
\end{lemma}

We can express $R^n$ by $\delta^n$, $\gamma_n$ and $d^n$ as follows.  

\begin{lemma} \label{R^n}
For any $n \geq 1$,
\[
R^n (l) = \dfrac{\gamma_n + l d^n}{\delta^n}.
\]
\end{lemma}


\begin{remark}
Let $f_0(z,w) = (z^{\delta}, z^{\gamma} w^d)$.
Then $f_0(z, cz^l) = (z^{\delta}, c^d z^{\gamma + l d})$
and 
\[
f_0^n (z, cz^l) = (z^{\delta^n}, c^{d^n} z^{\gamma_n + l d^n})
\]
since $f_0^n (z,w) = (z^{\delta^n}, z^{\gamma_n} w^{d^n})$.
Note that the ratio of the degrees with respect to $z$ of
the second and first components of $f_0^n (z, cz^l)$
coincides with $R^n(l)$.
In particular,
$f_0$ maps $\{ w = c z^l \}$ to $\{ w = c^d z^{R(l)} \}$,
and $f_0^n$ maps $\{ w = c z^l \}$ to $\{ w = c^{d^n} z^{R^n (l)} \}$.
\end{remark}

\subsection{Proof of Theorem \ref{main thm of weights when d>0}}

By definition,
$w(q) = \gamma + l d$ for any $l$ in $\mathcal{I}_f$.
Moreover,
using the previous lemmas,
we prove that
$w(Q^n) = \gamma_n + l d^n$ for any $l$ in $\mathcal{I}_f$
and for any $n \geq 1$.


\begin{proof}[Proof of Theorem \ref{main thm of weights when d>0}]
Fix any $l$ in $\mathcal{I}_f$
and let $w = w_l$.
For a monomial term $z^i w^j$,
we call $i+lj$ the weight of $z^i w^j$.
As in the description of the theorem,
we omit the coefficient of the term $z^{\gamma_n} w^{d^n}$
for simplicity.

We first show the equality $w(Q^2) = \gamma_2 + l d^2$.
Note that 
\[
Q^2(z,w) = \sum b_{ij} \big( p(z) \big)^i \big( q(z,w) \big)^j.
\]
Since $w(p) = \delta$ and $w(q) = \gamma + l d$,
$w(Q^2) 
\geq \min \{ \delta i + (\gamma + l d) j \, | \, b_{ij} \neq 0 \}$. 
It follows from Lemma \ref{dyn of R when d>0} 
that $i + R(l) j \geq \gamma + R(l) d$,
and so
$\delta i + (\gamma + l d) j \geq \delta \gamma + (\gamma + l d) d$ 
by Lemma \ref{R^n}.
Therefore,
$w(Q^2) \geq \delta \gamma + (\gamma + l d) d 
= (\delta + d) \gamma + l d^2 = \gamma_2 + l d^2$.
On the other hand,
the term 
\[
\big( z^{\delta} \big)^{\gamma} \big( z^{\gamma} w^d \big)^d = z^{\gamma_2} w^{d^2}
\]
has weight $\gamma_2 + l d^2$ and,
moreover,
has the smallest degree $d^2$ with respect to $w$ in $Q^2$. 
Because $d > 0$,
there are no other choices 
than $(z^{\delta})^{\gamma} (z^{\gamma} w^d)^d$
that generate terms of bidegree $(\gamma_2, d^2)$.
Therefore,
$Q^2$ contains the term $z^{\gamma_2} w^{d^2}$ and so
$w(Q^2) = w(z^{\gamma_2} w^{d^2}) = \gamma_2 + l d^2$.

We next show the equality $w(Q^3) = \gamma_3 + l d^3$ by the same strategy.
Note that 
\[
Q^3(z,w) = \sum b_{ij} \big( p^2(z) \big)^i \big( Q^2(z,w) \big)^j.
\]
Since $w(p^2) = \delta^2$ and $w(Q^2) = \gamma_2 + l d^2$,
$w(Q^3) \geq \min \{ \delta^2 i + (\gamma_2 + l d^2) j \, | \, b_{ij} \neq 0 \}$.
It follows from Lemma \ref{dyn of R when d>0}  
that $i + R^2(l) j \geq \gamma + R^2(l) d$,
and so
$\delta^2 i + (\gamma_2 + l d^2) j \geq \delta^2 \gamma + (\gamma_2 + l d^2) d$ 
by Lemma \ref{R^n}.
Therefore,
$w(Q^3) \geq \delta^2 \gamma + (\gamma_2 + l d^2) d 
= (\delta^2 \gamma + d \gamma_2) + l d^3 = \gamma_3 + l d^3$.
On the other hand,
the term 
\[
\big( z^{\delta^2} \big)^{\gamma} \big( z^{\gamma_2} w^{d^2} \big)^d = z^{\gamma_3} w^{d^3}
\]
has weight $\gamma_3 + l d^3$ and,
moreover,
has the smallest degree $d^3$ with respect to $w$ in $Q^3$. 
Because $d > 0$,
there are no other choices 
than $(z^{\delta^2})^{\gamma} (z^{\gamma_2} w^{d^2})^d$
that generate terms of bidegree $(\gamma_3, d^3)$.
Therefore,
$Q^3$ contains the term $z^{\gamma_3} w^{d^3}$ and so
$w(Q^3) = w(z^{\gamma_3} w^{d^3}) = \gamma_3 + l d^3$.

By repeating this process,
it follows that $w(Q^n) \geq \gamma_n + l d^n$ and
that $Q^n$ contains the term $z^{\gamma_n} w^{d^n}$
for any $n \geq 1$.
Therefore, 
$w(Q^n) = w(z^{\gamma_n} w^{d^n}) = \gamma_n + l d^n$.
\end{proof}

We remark that 
one can show that the coefficient of the term $z^{\gamma_n} w^{d^n}$ is
\[
a_{\delta}^{\gamma_{n-1} + \gamma_{n-2} + \cdots + \gamma_{2} + \gamma_{1}}
b_{\gamma d}^{d^{n-1} + d^{n-2} + \cdots + d + 1}
\]
from the construction of the term and 
the equality $\delta^k \gamma + d \gamma_k = \gamma_{k+1}$.

Because $w(Q^n)$ is the minimum $x$-intercept of the lines with slope $- l^{-1}$
that intersect $N(Q^n)$,
Theorem \ref{main thm of weights when d>0} implies that
$(\gamma_n, d^n)$ belongs to the boundary of $N(Q^n)$.
Moreover,
we obtain the following corollary
because the smallest degree with respect to $w$ of all terms in $Q^n$ is $d^n$.

\begin{corollary}
The bidegree
$(\gamma_n, d^n)$ is the vertex of $N(Q^n)$
whose $y$-coordinate is minimum for any $n \geq 1$. 
\end{corollary}

Therefore,
$N(Q^n)$ is included in the upper-right region
that is surrounded by the line with slope $-l_1^{-1}$
and the horizontal line,
which intersect at $(\gamma_n, d^n)$.

Let $- M_n$ be the slope of the line passing through 
$(\gamma_n, d^n)$ and the previous vertex if it exists,
let $M_n = \infty$ if $N(Q^n) = D(\gamma_n, d^n)$,
and let $M = M_1$.
Then $M = l_1^{-1}$,
and Theorem \ref{main thm of weights when d>0} implies the following corollary and
Corollary \ref{main cor on attr rates when d>0}  
because 
$w(Q^n)$ and $c(Q^n)$ are the minimum $x$-intercepts of 
the lines with slopes $-l^{-1}$ and $-1$ 
that intersect $N(Q^n)$. 

\begin{corollary}
Let $d > 0$. Then $M_n \geq M$ for any $n \geq 1$.
\end{corollary}

\subsection{Previous vertices and proof of Theorem \ref{main thm of attr rates when d>0}}

Furthermore,
we can show the existence of the vertex of $N(Q^n)$ 
that is previous to $(\gamma_n, d^n)$,
and specify it,
which induces the equality $M_n = M$.
Let $(A, B) = (n_{s-1}, m_{s-1})$,
\begin{align*}
(A_n, B_n) 
&= ((\delta^{n-1} + \delta^{n-2} d + \cdots + \delta d^{n-2}) \gamma + A d^{n-1}, B d^{n-1}) \\
&= (\gamma_n - (\gamma - A) d^{n-1}, B d^{n-1}) \text{ and} \\
(A_n^*, B_n^*) 
&= ( (\delta^{n-1} + \delta^{n-2} B + \cdots + \delta B^{n-2} + B^{n-1}) A,  B^{n}).
\end{align*}

\begin{proposition} \label{A_n and B_n}
Let $d> 0$.
If $\delta < T_{s-1}$, then
$(A_n, B_n)$ is the vertex of $N(Q^n)$ 
that is previous to $(\gamma_n, d^n)$,
$M_n = M$, and
$\delta^n$ is smaller than the y-intercept of 
the line passing through $(\gamma_n, d^n)$ and  $(A_n, B_n)$
for any $n \geq 1$.
If $\delta = T_{s-1}$, then
$(A_n^*, B_n^*)$ is the vertex of $N(Q^n)$ 
that is previous to $(\gamma_n, d^n)$,
$M_n = M$, and
$\delta^n$ coincides with the y-intercept of 
the line passing through $(\gamma_n, d^n)$ and  $(A_n^*, B_n^*)$
for any $n \geq 1$.
\end{proposition}

Before go into the proof,
we show an idea of how to find out $(A_2, B_2)$ and $(A_2^*, B_2^*)$.
Let us pick up a term $(z^{\delta})^I (z^i w^j)^J$ in $Q^2$, 
where $(i,j) \neq (\gamma, d)$.
As the same as the proof of Theorem \ref{main thm of weights when d>0},
\[
w( (z^{\delta})^I (z^i w^j)^J )
= \delta I + (i + l j) J
\geq \delta I + (\gamma + l d) J
\geq \delta \gamma + (\gamma + l d) d
= \gamma_2 + l d^2
\] 
for any $l$ in $\mathcal{I}_f$.
The equality in the first inequality holds 
if and only if $i + lj = \gamma + ld$ and $l = l_1$.
Let $l = l_1$ hereafter.
Note that the set 
\[
\{ (i,j) \, | \, i + l_1 j = \gamma + l_1 d \text{ and } b_{ij} \neq 0 \}
\]
lies on the side of $N(q)$ 
whose ends are $(n_{s-1}, m_{s-1})$ and $(\gamma, d)$. 
If $\delta < T_{s-1}$, then
the equality in the second inequality,
which is rewritten as $I + R(l_1) J = \gamma + R(l_1) d$, 
holds if and only if $(I,J) = (\gamma, d)$, 
since $R(l_1) > l_1$. 
Therefore,
\[
\big( z^{\delta} \big)^{\gamma} \big( z^{n_{s-1}} w^{m_{s-1}} \big)^d
\]
is the desired term in $Q^2$ and so
$(A_2, B_2) = (\delta \gamma + n_{s-1} d, m_{s-1} d)$.
On the other hand,
if $\delta = T_{s-1}$, then
the equality in the second inequality holds 
if and only if $I + l_1 J = \gamma + l_1 d$,
since $R(l_1) = l_1$. 
Therefore,
\[
\big( z^{\delta} \big)^{n_{s-1}} \big( z^{n_{s-1}} w^{m_{s-1}} \big)^{m_{s-1}}
\]
is the desired term in $Q^2$ and so
$(A_2^*, B_2^*) = (\delta n_{s-1} + n_{s-1} m_{s-1}, m_{s-1}^2)$.

Let $H^n_l$ be the polynomial that consists of all the terms in $Q^n$
of the smallest weight $\gamma_n + l d^n$. 
The explanation above indicates
the following proposition. 

\begin{proposition} \label{H^2 for case 2}
For any $l$ in $\mathcal{I}_f$,
$H_{l}^1 (z,w) = \sum_{i + l j= \gamma + l d} b_{ij} z^i w^j$ and
\[
H_{l}^2 (z,w) 
= \sum_{I + R(l) J = \gamma + R(l) d} b_{IJ}
\left( a_{\delta} z^{\delta} \right)^{I} \bigg( \sum_{i + l j= \gamma + l d} b_{ij} z^i w^j \bigg)^{J}.
\]
More precisely,
if $\delta < T_{s-1}$,
then $\mathcal{I}_f = [l_1, \alpha]$ or $\mathcal{I}_f = [l_1, \infty)$,
\begin{align*}
H^2_{l_1} (z,w) 
&= b_{\gamma d} \left( a_{\delta} z^{\delta} \right)^{\gamma} 
           \bigg( \sum_{i+l_1 j= \gamma +l_1 d} b_{ij} z^i w^j \bigg)^d \\
&= b_{\gamma d} \left( a_{\delta} z^{\delta} \right)^{\gamma} 
  \left( b_{\gamma d} z^{\gamma} w^d + \cdots + b_{AB} z^A w^B \right)^d
  \text{ and} \\ 
H_{l}^2 (z,w) 
&= b_{\gamma d} \left( a_{\delta} z^{\delta} \right)^{\gamma} \left( b_{\gamma d} z^{\gamma} w^d \right)^d
\text{ for any $l$ in } \mathcal{I}_f \setminus \{ l_1 \}.
\end{align*}
On the other hand,
if $\delta = T_{s-1}$,
then $\mathcal{I}_f = \{ l_1 \}$ and
\begin{align*}
H_{l_1}^2 (z,w) 
&= \sum_{I+l_1J= \gamma +l_1 d} b_{IJ}
 \left( a_{\delta} z^{\delta} \right)^{I} \bigg( \sum_{i+l_1 j= \gamma +l_1 d} b_{ij} z^i w^j \bigg)^{J} \\
&= b_{\gamma d} \left( a_{\delta} z^{\delta} \right)^{\gamma} 
  \left( b_{\gamma d} z^{\gamma} w^d + \cdots + b_{AB} z^A w^B \right)^d + \cdots \\
& \ \ \ \ + b_{AB} \left( a_{\delta} z^{\delta} \right)^{A} 
  \left( b_{\gamma d} z^{\gamma} w^d + \cdots + b_{AB} z^A w^B \right)^B.
\end{align*}
\end{proposition}


The lemma below follows immediately from the definition of $\gamma_n$.

\begin{lemma} \label{slope lem}
The slopes of the lines passing through $(0, \delta^n)$ and $(\gamma_n, d^n)$
for all $n \geq 1$ are all the same: 
$(\delta^n - d^n)/\gamma_n = (\delta - d)/\gamma$.
\end{lemma}

As a consequence of Proposition \ref{A_n and B_n} 
and Lemma \ref{slope lem},
it follows that
\begin{align*}
\dfrac{B_n - d^n}{\gamma_n - A_n} &= \dfrac{B - d}{\gamma - A} 
> \dfrac{\delta - d}{\gamma} = \dfrac{\delta^n - d^n}{\gamma_n}
\text{ if $\delta < T_{s-1}$, and} \\
\dfrac{B_n^* - d^n}{\gamma_n - A_n^*} &= \dfrac{B - d}{\gamma - A} 
= \dfrac{\delta - d}{\gamma} = \dfrac{\delta^n - d^n}{\gamma_n} 
\text{ if $\delta = T_{s-1}$}.
\end{align*}

Now we are ready to prove Proposition \ref{A_n and B_n}.

\begin{proof}[Proof of Proposition \ref{A_n and B_n}]
As we saw in the explanation above, 
we have to choose $l_1$ as $l$ in $\mathcal{I}_f$
to find out the vertex of $N(Q^n)$ 
that is previous to $(\gamma_n, d^n)$.
Then $w_{l_1}(Q^n) = \gamma_n + l_1 d^n$
by Theorem \ref{main thm of weights when d>0}.

If $\delta < T_{s-1}$,
then it follows from Proposition \ref{H^2 for case 2} that the term
\[
\big( z^{\delta} \big)^{\gamma} \big( z^A w^B \big)^d 
= z^{\delta \gamma + Ad} w^{Bd}
\]
has weight $\gamma_2 + l_1 d^2$ and, moreover,
has the biggest degree with respect to $w$ in $H^2_{l_1}$.
Let 
\[
(A_2, B_2) = (\delta \gamma + Ad, Bd).
\]
Since there are no other choices 
that generate terms of bidegree $(A_2, B_2)$,
$Q^2$ contains the term $z^{A_2} w^{B_2}$ and so
$(A_2, B_2)$ is the vertex of $N(Q^2)$ 
that is previous to $(\gamma_2, d^2)$.
Since the weight of $z^{\gamma_2} w^{d^2}$ and $z^{A_2} w^{B_2}$ are the same,
$M_2 = M$.
Consequently,
\[
\dfrac{B_2 - d^2}{\gamma_2 - A_2} =  \dfrac{B - d}{\gamma - A}
> \dfrac{\delta - d}{\gamma} = \dfrac{\delta^2 - d^2}{\gamma_2}
\]
by Lemma \ref{slope lem}
and so $\delta^2$ is smaller than the $y$-intercept of 
the line passing through $(\gamma_2, d^2)$ and  $(A_2, B_2)$.

It follows from the same arguments that the term 
\[
\big( z^{\delta^2} \big)^{\gamma} \big( z^{A_2} w^{B_2} \big)^d = z^{\delta^2 \gamma + A_2 d} w^{B_2 d}
\] 
has weight $\gamma_3 + l_1 d^3$ and, moreover,
has the biggest degree with respect to $w$ in $H^3_{l_1}$.
Let
\[
(A_3, B_3) = (\delta^2 \gamma + \delta d \gamma + A d^2, B d^2).
\]
Then it follows from the same arguments
that $(A_3, B_3)$ is the vertex of $N(Q^3)$ that is previous to $(\gamma_3, d^3)$, 
$M_3 = M$, and
$\delta^3$ is smaller than the $y$-intercept of 
the line passing through $(0, \delta^3)$ and $(\gamma_3, d^3)$.
Repeating this process,
we obtain the required vertices and properties.

If $\delta = T_{s-1}$,
then it follows Proposition \ref{H^2 for case 2} that the term 
\[
\big( z^{\delta} \big)^{A} \big( z^A w^B \big)^{B}
= z^{\delta A + A B} w^{B^2}
\]
has weight $\gamma_2 + l_1 d^2$ and, moreover,
has the biggest degree with respect to $w$ in $H^2_{l_1}$.
Let
\[
(A_2^*, B_2^*) = (\delta A + A B, B^2) = ((\delta + B)A, B^2).
\]
Since there are no other choices 
that generate terms of bidegree $(A_2^*, B_2^*)$,
it is the vertex of $N(Q^2)$ 
that is previous to $(\gamma_2, d^2)$.
Moreover,
$M_2 = M$
and so 
$\delta^2$ coincides with the $y$-intercept of 
the line passing through $(\gamma_2, d^2)$ and $(A_2^*, B_2^*)$.

It follows from the same arguments that the term  
\[
\big( z^{\delta^2} \big)^{A} \big( z^{A_2^*} w^{B_2^*} \big)^B 
= z^{\delta^2 A + A_2^* B} w^{B_2^* B}
\] 
has weight $\gamma_3 + l_1 d^3$ and, moreover,
has the biggest degree with respect to $w$ in $H_{l_1}^3$.
Let
\[
(A_3^*, B_3^*) = (\delta^2 A + (\delta + B)AB, B^3) = ((\delta^2 + \delta B + B^2)A, B^3).
\]
Then it follows from the same arguments
that $(A_3^*, B_3^*)$ is the vertex of $N(Q^3)$ that is previous to $(\gamma_3, d^3)$,
$M_3 = M$, and
$\delta^3$ coincides with the $y$-intercept of 
the line passing through $(0, \delta^3)$ and $(\gamma_3, d^3)$.
Repeating this process,
we obtain the required vertices and properties.
\end{proof}


Theorem \ref{main thm of weights when d>0} and
Proposition \ref{A_n and B_n} induce Theorem \ref{main thm of attr rates when d>0}.

\begin{proof}[Proof of Theorem \ref{main thm of attr rates when d>0}]
We only show the improved parts.
Recall that 
$N(Q^n)$ is included in the upper-right region
that is surrounded by the line with slope $-l_1^{-1}$
and the horizontal line,
which intersect at $(\gamma_n, d^n)$,
and that
$c(Q^n)$ is the minimum $x$-intercept of 
the lines with slope $-1$ that intersect $N(Q^n)$. \\
(2)
If $\delta < T_{s-1}$,
then $c(Q^n) \leq A_n + B_n < \gamma_n + d^n$.
If $\delta = T_{s-1}$,
then $c(Q^n) \leq A_n^* + B_n^* < \gamma_n + d^n$. \\
(3)
If $n_1 > 0$ or $s > 2$, 
then $A = n_{s-1} > 0$.
Hence $A_n > 0$ or $A_n^* > 0$ 
since $\delta$, $\gamma$, $d$, $A$ and $B$ are all positive. 
Therefore,
$N(Q^n)$ does not contain the $y$-intercept $(0, l_1^{-1} \gamma_n + d^n)$
of the line passing through $(\gamma_n, d^n)$ with slope $-l_1^{-1}$,
and so 
$c(Q^n) > l_1^{-1} \gamma_n + d^n$. \\
(4) 
If $\delta < T_{s-1}$,
then $A_1 = A = n_1 = 0$ and $A_n > 0$ for any $n \geq 2$.
Hence $c(q) = B_1 = B = l_1^{-1} \gamma + d$, and
$c(Q^n) > l_1^{-1} \gamma_n + d^n$ for any $n \geq 2$
because $N(Q^n)$ does not contain $(0, l_1^{-1} \gamma_n + d^n)$. \\
(5)
If $\delta = T_{s-1}$,
then $A_n^* = 0$ and 
$c(Q^n) = B_n^* = B^n = \delta^n = l_1^{-1} \gamma_n + d^n$.
\end{proof}

\begin{remark}
If $d > 0$, 
then $\gamma_n \to \infty$, and $A_n$ or $A_n^* \to \infty$
as $n \to \infty$
since $\delta > 0$ and $\gamma > 0$.
\end{remark}

\begin{remark}
Recall that $f$ preserves the open wedge $U^l$ for any $l$ in $\mathcal{I}_f$.
Hence $f^n$ preserves $U^l$ for any $l$ in $\mathcal{I}_f$.
Therefore,
one may expect that $\mathcal{I}_f \subset \mathcal{I}_{f^n}$,
where $\mathcal{I}_{f^n}$ is the interval of the weights for $f^n$.
In fact, 
Proposition \ref{A_n and B_n} and Lemma \ref{slope lem} imply that
the equality $\mathcal{I}_f = \mathcal{I}_{f^n}$ holds
for any $n \geq 1$.
\end{remark}

\section{Shape of Newton polygon of $Q^n$ when $d = 0$}

We complete the investigation of $w(Q^n)$ and $c(Q^n)$ for Case 2 in this section,
assuming that $d = 0$.
We first give a rough explanation of 
the differences between the cases $d > 0$ and $d = 0$
and a summary of the results in Section 4.1. 
Some of the results are illustrated 
in terms of the blow-ups in Section 4.2. 
We then give a preliminary lemma in Section 4.3 and
more detailed explanations of the results
for the cases $\delta < T_{s-1}$ and $\delta = T_{s-1}$
in Sections 4.4 and 4.5, respectively,
although 
we omit the precise proofs of the results
because they are more or less similar to the case $d > 0$.

\subsection{Explanation of differences and Summary of results}

The situation for the case $d = 0$ is different from that for the previous case $d > 0$.
Let us give here a rough explanation.
If  $d = 0$,
then $(\gamma_n, d^n) = (\delta^{n-1} \gamma, 0)$.
Recall that, for the case $d > 0$,
the bidegree $(A_n, B_n)$ or $(A_n^*, B_n^*)$ is the vertex of $N(Q^n)$
that is previous to $(\gamma_n, d^n)$
if $\delta < T_{s-1}$ or $\delta = T_{s-1}$. 
If $d = 0$ and $\delta < T_{s-1}$,
then $(A_n, B_n)$ coincides with $(\gamma_n, 0)$,
and the equality $M_n = M$ does not hold;
in fact, the inequality $M_n > M$ holds for any $n \geq 2$.
However,
the term $z^{\gamma_n}$ remains forever, and 
the equality $w (Q^n) = w (z^{\gamma_n}) = \gamma_n$ holds
for any $n \geq 1$ and
for any $l$ in $\mathcal{I}_f$,
where $\mathcal{I}_f = [l_1, \alpha]$.
On the other hand,
if $d = 0$ and $\delta = T_{s-1}$,
then $(A_n, B_n)$ coincides with $(\gamma_n, 0)$ and, moreover,
the term $z^{\gamma_n}$ may vanish.
However,
the term $z^{A_n^*} w^{B_n^*}$ remains forever, and
the equality $w_{l_1} (Q^n) = w_{l_1} (z^{A_n^*} w^{B_n^*}) = \gamma_n$ holds
for any $n \geq 1$,
where $\mathcal{I}_f = \{ l_1 \}$.
Consequently,
even if $d = 0$,
we have the same equality on $w(Q^n)$ as the case $d > 0$.

\begin{theorem} \label{main thm of weights when d=0}
Let $d = 0$.
Then $w(Q^n) = \gamma_n$
for any $n \geq 1$ and for any $l$ in $\mathcal{I}_f$,
where $\mathcal{I}_f = [l_1, \alpha]$.
\end{theorem}

Unlike the case $d>0$,
$c(Q^n)$ can be bigger than $\gamma_n$
since the term $z^{\gamma_n}$ may vanish if $d = 0$.

\begin{corollary} \label{main cor on attr rates when d=0}
Let $d = 0$.
Then it follows for any $n \geq 1$ that 
\begin{enumerate}
\item[$(1)$] $\gamma_n \leq c(Q^n) \leq l_1^{-1} \gamma_n$
if $l_1 \leq 1$, and 
\item[$(2)$] $l_1^{-1} \gamma_n \leq c(Q^n) \leq \gamma_n$
if $l_1 > 1$.
\end{enumerate}
\end{corollary}

Moreover,
we have the following improved estimates on $c(Q^n)$.

\begin{theorem} \label{main thm of attr rates when d=0 and delta < T}
Let $d = 0$ and $\delta < T_{s-1}$.
Then
\begin{enumerate}
\item[$(1)$] $c(Q^n) = \gamma_n$
for any $n \geq 1$ if $l_1 \leq 1$, and 
\item[$(2)$] $l_1^{-1} \gamma \leq c(q)  < \gamma$ \text{ and } 
        $l_1^{-1} \gamma_n < c(Q^n) \leq \gamma_n$ for any $n \geq 2$ if $l_1 > 1$.
\end{enumerate}
\end{theorem}

\begin{theorem} \label{main thm of attr rates when d=0 and delta = T}
Let $d = 0$ and $\delta = T_{s-1}$.
Then it follows for any $n \geq 1$ that
\begin{enumerate}
\item[$(1)$] $\gamma_n \leq c(Q^n)  \leq l_1^{-1} \gamma_n$
if $l_1 \leq 1$, and 
\item[$(2)$] $l_1^{-1} \gamma_n \leq c(Q^n)  < \gamma_n$
if $l_1 > 1$.  
\end{enumerate}
More precisely,  if $l_1 < 1$ and the term $z^{\gamma_n}$ vanishes for some $n_0$, then
\begin{enumerate}
\item[$(3)$] $\gamma_n < c(Q^n)$ for any $n \geq n_0$.
\end{enumerate}
\end{theorem}

\subsection{Partial illustration of results in terms of blow-ups}

Assuming that $l_1$ is an integer,
we illustrate Theorem \ref{main thm of weights when d=0} and
Corollary \ref{main cor on attr rates when d=0} partially in terms of the blow-up. 
Let $\pi_1 (z,c) = (z, z^{l_1} c)$
and $\tilde{f} = \pi_1^{-1} \circ f \circ \pi_1$.
By Lemma \ref{blow-up lem for case 2},
$N(\tilde{q}) = D(\tilde{\gamma}, 0)$.
If $\delta < T_{s-1}$,
then $\tilde{\gamma} > 0$ and the same claims hold as the case $d > 0$:
$\tilde{Q}^n$ contains the term $z^{\tilde{\gamma}_n}$ and
$N(\tilde{Q}^n) = D(\tilde{\gamma}_n, 0)$,
which imply 
Theorem \ref{main thm of weights when d=0}
and Corollary \ref{main cor on attr rates when d=0}.
On the other hand,
if $\delta = T_{s-1}$,
then $\tilde{\gamma} = 0$ and 
the term $z^{\tilde{\gamma}_n}$ may vanish.
However,
the inclusion $N(\tilde{Q}^n) \subset D(\tilde{\gamma}_n, 0)$ still holds.
Hence
$N(Q^n)$ is included in the upper-right region
that is surrounded by the line with slope $-l_1^{-1}$
and the $x$-axis,
which intersect at $(\gamma_n, 0)$. 
Therefore,
we obtain the estimates on $w(Q^n)$ and $c(Q^n)$
from the below: 
$w(Q^n) \geq \gamma_n$,
$c(Q^n) \geq \gamma_n$ if $l_1 \leq 1$, and
$c(Q^n) \geq l_1^{-1} \gamma_n$ if $l_1 > 1$
for any $n \geq 1$ and for any $l$ in $\mathcal{I}_f$.


\subsection{Preliminary lemma: Affine dynamics on the interval}

If $d = 0$,
then $R \equiv \alpha$.
Hence $R$ collapses any point to $\alpha$, and
$\mathcal{I}_f = [l_1, \alpha]$
since $\delta > d$.

\begin{lemma}
For any $l$ in $\mathcal{I}_f$,
it follows that $R^n (l)$ belongs to $\mathcal{I}_f$ 
for any $n \geq 1$.
More precisely,
$R \equiv \alpha$.
\end{lemma}

\subsection{The case $d = 0$ and $\delta < T_{s-1}$}

We showed in the proof of Theorem \ref{main thm of attr rates when d>0} that, 
if $d > 0$ and $\delta < T_{s-1}$,
then 
the polynomial $H^2_{l_1}$, 
\[
b_{\gamma d} \big( a_{\delta} z^{\delta} \big)^{\gamma} 
  \big( b_{\gamma d} z^{\gamma} w^d + \cdots + b_{AB} z^A w^B \big)^d,
\]
generates the important terms $z^{\gamma_2} w^{d^2}$ and $z^{A_2} w^{B_2}$.
However,
if $d = 0$ and $\delta < T_{s-1}$,
then these terms coincide:
\[
b_{\gamma 0} \big( a_{\delta} z^{\delta} \big)^{\gamma} 
  \big( b_{\gamma 0} z^{\gamma} w^0 + \cdots + b_{AB} z^A w^B \big)^0
= a_{\delta}^{\gamma} b_{\gamma 0} z^{\gamma_2}. 
\]
On the other hand, 
it follows from the same arguments as the case $d > 0$ that,
if $d = 0$ and $\delta < T_{s-1}$,
then the term $z^{\gamma_n}$ remains forever
as the unique term of the smallest weight $\gamma_n$. 
Hence we obtain the following three propositions,
which imply Theorem \ref{main thm of weights when d=0} 
for the case $\delta < T_{s-1}$
and Theorem \ref{main thm of attr rates when d=0 and delta < T}.

\begin{proposition}
Let $d = 0$ and $\delta < T_{s-1}$.
Then $Q^n$ contains the term $z^{\gamma_n}$ and 
$w(Q^n) = \gamma_n$ 
for any $n \geq 1$ and 
for any $l$ in $\mathcal{I}_f$, where $\mathcal{I}_f = [l_1, \alpha]$.
\end{proposition}

\begin{proposition}
Let $d = 0$ and $\delta < T_{s-1}$.
Then $M_n > M$ for any $n \geq 2$.
\end{proposition}

\begin{proposition}
Let $d = 0$ and $\delta < T_{s-1}$. Then
\begin{enumerate}
\item[$(1)$] $c(Q^n) = \gamma_n$ for any $n \geq 1$ if $l_1 \leq 1$, and 
\item[$(2)$] $l_1^{-1} \gamma \leq c(q) \leq A + B < \gamma$ \text{ and } 
        $l_1^{-1} \gamma_n < c(Q^n) \leq \gamma_n$ for any $n \geq 2$ if $l_1 > 1$.
\end{enumerate}
\end{proposition}

\begin{remark}[Unboundedness of $M_n$]
Let $d = 0$ and $\delta < T_{s-1}$.
Then
\[
M < M_n \leq \dfrac{B_n^*}{\gamma_n - A_n^*}
\]
for any $n \geq 2$.
If $A = 0$, then $(A_n^*, B_n^*) = (0, B^n)$.
Hence $B_n^*/\gamma_n \to \infty$ as $n \to \infty$.
Moreover,
if $(A_n^*, B_n^*)$ is the vertex of $N(Q^n)$ that is previous to $(\gamma_n, d^n)$,
then $M_n = B_n^*/\gamma_n \to \infty$ as $n \to \infty$.
On the other hand,
if $A > 0$, then $A_n^* > \gamma_n$ for any large $n$.
More strongly, it can happen that
$M_n = \infty$ for any large $n$.
\end{remark}

\subsection{The case $d = 0$ and $\delta = T_{s-1}$}

We showed in the proof of Theorem \ref{main thm of attr rates when d>0} that, 
if $d > 0$ and $\delta = T_{s-1}$,
then 
the polynomial $H^2_{l_1}$,
\begin{align*}
& b_{\gamma d} \big( a_{\delta} z^{\delta} \big)^{\gamma} 
  \big( b_{\gamma d} z^{\gamma} w^d + \cdots + b_{AB} z^A w^B \big)^d + \cdots \\
& + b_{AB} \big( a_{\delta} z^{\delta} \big)^{A} 
  \big( b_{\gamma d} z^{\gamma} w^d + \cdots + b_{AB} z^A w^B \big)^B,
\end{align*}
generates the important terms $z^{\gamma_2} w^{d^2}$ and $z^{A_2^*} w^{B_2^*}$.
However,
if $d = 0$ and $\delta = T_{s-1}$,
then $H^2_{l_1}$ has the terms that consist only of $z$
other than $a_{\delta}^{\gamma} b_{\gamma 0} z^{\gamma_2}$:
\begin{align*}
& b_{\gamma 0} \big( a_{\delta} z^{\delta} \big)^{\gamma} 
  \big( b_{\gamma 0} z^{\gamma} w^0 + \cdots + b_{AB} z^A w^B \big)^0 + \cdots \\
& + b_{AB} \big( a_{\delta} z^{\delta} \big)^{A} 
  \big( b_{\gamma 0} z^{\gamma} w^0 + \cdots + b_{AB} z^A w^B \big)^B \\
= \ & b_{\gamma 0} \big( a_{\delta} z^{\delta} \big)^{\gamma} + \cdots 
+ b_{AB} \big( a_{\delta} z^{\delta} \big)^{A} 
  \Big\{ \big( b_{\gamma 0} z^{\gamma} \big)^B + \cdots + \big( b_{AB} z^A w^B \big)^B \Big\} \\
= \ & a_{\delta}^{\gamma} b_{\gamma 0} z^{\delta \gamma} + \cdots 
 + a_{\delta}^A b_{AB} b_{\gamma 0}^B z^{\delta A + \gamma B} + \cdots 
 + a_{\delta}^A b_{AB}^{B+1} z^{\delta A + AB} w^{B^2}.
\end{align*}
Note that $H^2_{l_1}$ contains the term 
$b_{IJ} (a_{\delta} z^{\delta} )^I (b_{\gamma 0} z^{\gamma})^J
= a_{\delta}^I b_{IJ} b_{\gamma 0}^J z^{\delta I + \gamma J}$
for any $(I,J)$ such that $I + l_1J= \gamma + l_1 d$.
Because the weights of all the terms in $H^2_{l_1}$ are the same $\gamma_2$,
it follows that
\[
\delta \gamma = \delta A + \gamma B = \delta I + \gamma J.
\]
We remark that this equality also follows directly from the condition $\delta = T_{s-1}$
since
$T_{s-1} = \gamma B/(\gamma - A) = \gamma J/(\gamma - I)$.
Hence
\[
H^2_{l_1} (z,w) = \bigg\{ \sum_{I+l_1J= \gamma} a_{\delta}^I b_{IJ} b_{\gamma 0}^J \bigg\} 
z^{\delta \gamma} + \cdots + a_{\delta}^A b_{AB}^{B+1} z^{A_2^*} w^{B_2^*}. 
\]
Therefore,
the term $z^{\gamma_2}$ vanishes if 
\[
\sum_{I+l_1J= \gamma} a_{\delta}^I b_{IJ} b_{\gamma 0}^J = 0. 
\]
If the term $z^{\gamma_n}$ remains forever,
then the same equality on $w(Q^n)$ and inequalities on $c(Q^n)$ hold
as the case $d > 0$.
Note that,
for example,
if $a_{\delta}$ and
all the non-zero coefficients of $q$ are positive,
then the term $z^{\gamma_n}$ remains forever.
Even if $z^{\gamma_n}$ vanishes for some $n = n_0$,
the term $z^{A_n^*} w^{B_n^*}$ remains forever and
$w_{l_1} (Q^n) = w_{l_1} (z^{A_n^*} w^{B_n^*}) = \gamma_n$
for any $n \geq n_0$.
In particular,
we obtain the following two propositions, 
which imply Theorem \ref{main thm of weights when d=0} 
for the case $\delta = T_{s-1}$
and Theorem \ref{main thm of attr rates when d=0 and delta = T}.

\begin{proposition}
Let $d = 0$ and $\delta = T_{s-1}$.
Then $w_{l_1} (Q^n) = \gamma_n$ for any $n \geq 1$ and $\mathcal{I}_f = \{ l_1 \}$.
\end{proposition}

\begin{proposition}
Let $d = 0$ and $\delta = T_{s-1}$.
Then it follows for any $n \geq 1$ that
\begin{enumerate}
\item[$(1)$] $\gamma_n \leq c(Q^n) \leq A_n^* + B_n^* \leq l_1^{-1} \gamma_n$
if $l_1 \leq 1$, and 
\item[$(2)$] $l_1^{-1} \gamma_n \leq c(Q^n) \leq A_n^* + B_n^* < \gamma_n$
if $l_1 > 1$.
\end{enumerate}
More precisely,  if $l_1 < 1$ and the term $z^{\gamma_n}$ vanishes for some $n_0$, then
\begin{enumerate}
\item[$(3)$] $\gamma_n < c(Q^n)$ for any $n \geq n_0$.
\end{enumerate}
\end{proposition}

\section{Shape of Newton polygon of $Q^n$ for Case 3}

In this section we deal with Case 3.
Let $s > 1$,
\[
T_1 \leq \delta, \
(\gamma, d) = (n_1, m_1) \text{ and }
l_2 = \frac{n_2 - n_1}{m_1 - m_2}. 
\]
Note that $\delta \geq d > 0$ by the setting.
We first give a summary of the results in Section 5.1
and illustrate some of them in terms of the blow-ups in Section 5.2. 
Preparing a lemma in Section 5.3,
we prove Theorem \ref{main thm of weights for case 3} below in Section 5.4.
Finally,
we show the existence of the vertex of $N(Q^n)$
that is next to $(\gamma_n, d^n)$
in most cases
and specify it in Section 5.5,
which induces Theorem \ref{main thm of attr rates for case 3} below.

\subsection{Summary of results}

\begin{theorem} \label{main thm of weights for case 3}
It follows that 
$Q^n$ contains the term $z^{\gamma_n} w^{d^n}$ and
$w(Q^n) = \gamma_n + l d^n$ 
for any $n \geq 1$ and
for any $l$ in $\mathcal{I}_f$,
where
$\mathcal{I}_f = [\alpha, l_2]$ or $\mathcal{I}_f = (0, l_2]$
if $\gamma > 0$ or $\gamma = 0$.
\end{theorem}

Moreover,
it follows that
$(\gamma_n, d^n)$ is the vertex of $N(Q^n)$
whose $x$-coordinate is minimum.
Hence the theorem above implies the following corollary. 

\begin{corollary} \label{main cor on attr rates for case 3}
It follows for any $n \geq 1$ that
\begin{enumerate}
\item[$(1)$] $c(Q^n) = \gamma_n + d^n$ 
if $l_2 \geq 1$, and
\item[$(2)$] $\gamma_n + l_2 d^n \leq c(Q^n) \leq \gamma_n + d^n$
if $l_2 < 1$.
\end{enumerate}
\end{corollary}

Let $(C, D) = (n_2, m_2)$ and
we define $(C_n, D_n)$ and $(C_n^*, D_n^*)$ as the same as Case 2.
If $\delta > T_1$,
then the term $z^{C_n} w^{D_n}$ remains forever and
$(C_n, D_n)$ is the vertex of $N(Q^n)$ 
that is next to $(\gamma_n, d^n)$.
Similarly,
if $\delta = T_1$ and $m_2 > 0$,
then the term $z^{C_n^*} w^{D_n^*}$ 
remains forever and
$(C_n^*, D_n^*)$ is the vertex of $N(Q^n)$
that is next to $(\gamma_n, d^n)$.
On the other hand,
if $\delta = T_1$ and $m_2 = 0$,
then the term $z^{C_n^*} w^{D_n^*}$ may vanish.
By using these vertices,
we can improve the inequalities for the case $l_1 < 1$
in the corollary above
as follows.

\begin{theorem} \label{main thm of attr rates for case 3}
It follows for any $n \geq 1$ that
\begin{enumerate}
\item[$(1)$] $c(Q^n) = \gamma_n + d^n$
if $l_2 \geq 1$, and
\item[$(2)$] $\gamma_n + l_2 d^n \leq c(Q^n) \leq \gamma_n + d^n$
if $l_2 < 1$.
\end{enumerate}
Let $l_2 < 1$ and $m_2 > 0$.
Then it follows for any $n \geq 1$ that
\begin{enumerate}
\item[$(3)$] $\gamma_n + l_2 d^n < c(Q^n) < \gamma_n + d^n$.
\end{enumerate}
Let $l_2 < 1$ and $m_2 = 0$.
Then $s = 2$ and
\begin{enumerate}
\item[$(4)$] $\gamma_n + l_2 d^n = c(Q^n)$
for any $n \geq 1$ if $\delta > T_1$, 
\item[$(5)$] $\gamma_n + l_2 d^n = c(Q^n)$
for any $n \geq 1$ if $\delta = T_1$ and 
the term $z^{C_n^*} w^{D_n^*}$ remains forever, and
\item[$(6)$] $\gamma_n + l_2 d^n < c(Q^n)$ for any $n \geq n_0$
        if $\delta = T_1$ and 
the term $z^{C_n^*} w^{D_n^*}$ vanishes for some $n_0$. 
\end{enumerate}
\end{theorem}

\subsection{Illustration of results in terms of blow-ups}

Assuming that $l_2^{-1}$ is an integer,
we illustrate Theorem \ref{main thm of weights for case 3} and 
Corollary \ref{main cor on attr rates for case 3} in terms of the blow-up. 
Let $\pi_2 (t,w) = (t w^{l_2^{-1}}, w)$ 
and $\tilde{f} = \pi_2^{-1} \circ f \circ \pi_2$
as in Section 2.2.
By Lemma \ref{blow-up lem for case 3},
$N(\tilde{q}) = D(\gamma, \tilde{d})$,
where $\tilde{d} = l_2^{-1} \gamma + d$.
Although $\tilde{f}$ may not be skew product,
it is close to a monomial map 
by Proposition \ref{blow-up prop for case 3} 
and so $N(\tilde{Q}^n)$ has the unique vertex,
that should be $(\gamma_n, \tilde{d}^n)$.
Hence 
$N(\tilde{Q}^n) = D(\gamma_n, \tilde{d}^n)$.
This implies that 
$Q^n$ contains the term $z^{\gamma_n} w^{d^n}$, and that
$N(Q^n)$ is included in the upper-right region
that is surrounded by the line with slope $-l_2^{-1}$
and the vertical line,
which intersect at $(\gamma_n, d^n)$. 
Therefore,
we obtain Theorem \ref{main thm of weights for case 3}
and Corollary \ref{main cor on attr rates for case 3} when $l_2^{-1}$ is an integer,
because
$w(Q^n)$ and $c(Q^n)$ are the minimum $x$-intercepts of 
the lines with slopes $-l^{-1}$ and $-1$ 
that intersect $N(Q^n)$. 

\subsection{Preliminary lemma: Affine dynamics on the interval}

If $\delta > d$,
then $\gamma \geq 0$ and
$R$ is a contracting function around the fixed point $\alpha$,
where $\alpha \geq 0$.
If $\delta = d$,
then $\gamma = 0$ and so $R(l) = l$.
Since $\mathcal{I}_f = [\alpha, l_2]$ or $\mathcal{I}_f = (0, l_2]$
if $\gamma > 0$ or $\gamma = 0$, 
we can conclude as follows.

\begin{lemma} \label{dyn of R for case 3}
For any $l$ in $\mathcal{I}_f$,
it follows that $R^n (l)$ belongs to $\mathcal{I}_f$ for any $n \geq 1$ and
the sequence $\{ R^n (l) \}_{n \geq 1}$ is decreasing or $R = id$. 
More precisely, 
$R^n (l) \to \alpha$ as $n \to \infty$ or $R = id$ 
if $\delta > d$ or $\delta = d$. 
\end{lemma}

\subsection{Proof of Theorem \ref{main thm of weights for case 3}}

We can prove Theorem \ref{main thm of weights for case 3}  
by combining arguments similar to the proofs of
Theorem \ref{main thm of weights when d>0} and 
Proposition \ref{A_n and B_n} for Case 2.
The inequality $w(Q^n) \geq \gamma_n + l d^n$ 
for any $n \geq 1$ and any $l$ in $\mathcal{I}_f$ 
follows from Lemma \ref{dyn of R for case 3} and
the same arguments as the proof of
Theorem \ref{main thm of weights when d>0}. 
Although $d^n$ is not the smallest nor biggest degree
with respect to $w$ in $Q^n$,
it is the biggest degree with respect to $w$ in $H^n_l$;
similar arguments can be found in the proof of
Proposition \ref{A_n and B_n}.
We obtain the following proposition
by arguments similar to Case 2.


\begin{proposition} \label{H^2 for case 3}
If $\delta > T_{1}$,
then $\mathcal{I}_f = [ \alpha, l_2 ]$ or $\mathcal{I}_f = (0, l_2]$,
\begin{align*}
H^2_{l_2} (z,w) &= b_{\gamma d} \left( a_{\delta} z^{\delta} \right)^{\gamma} 
           \bigg( \sum_{i + l_2 j= \gamma + l_2 d} b_{ij} z^i w^j \bigg)^d \\
&= b_{\gamma d} \big( a_{\delta} z^{\delta} \big)^{\gamma} 
  \big( b_{\gamma d} z^{\gamma} w^d + \cdots + b_{CD} z^C w^D \big)^d
  \text{ and} \\
H^2_{l} (z,w) 
&= b_{\gamma d} \big( a_{\delta} z^{\delta} \big)^{\gamma} \big( b_{\gamma d} z^{\gamma} w^d \big)^d
\text{ for any $l$ in } \mathcal{I}_f - \{ l_2 \}.
\end{align*}
If $\delta = T_{1}$,
then $\mathcal{I}_f = \{ l_2 \}$ and
\begin{align*}
H^2_{l} (z,w) &= \sum_{I+l_2J= \gamma +l_2 d} b_{IJ}
 \left( a_{\delta} z^{\delta} \right)^{I} \bigg( \sum_{i+l_2 j= \gamma + l_2 d} b_{ij} z^i w^j \bigg)^{J} \\
&= b_{\gamma d} \big( a_{\delta} z^{\delta} \big)^{\gamma} 
  \big( b_{\gamma d} z^{\gamma} w^d + \cdots + b_{CD} z^C w^D \big)^d + \cdots \\
& \ \ \ \ + b_{CD} \big( a_{\delta} z^{\delta} \big)^{C} 
  \big( b_{\gamma d} z^{\gamma} w^d + \cdots + b_{CD} z^C w^D \big)^D.
\end{align*}
\end{proposition}

\begin{proof}[Proof of Theorem \ref{main thm of weights for case 3}]
Fix any $l$ in $\mathcal{I}_f$.
As mentioned above,
we have the inequality $w(Q^n) \geq \gamma_n + l d^n$ for any $n \geq 1$.

It follows from  Proposition \ref{H^2 for case 3} that
$z^{\gamma_2} w^{d^2}$
is the term in $H^2_{l}$
with the biggest degree $d^2$ with respect to $w$ in $H^2_{l}$;
hence $Q^2$ contains the term $z^{\gamma_2} w^{d^2}$.
It follows from the same arguments that 
$z^{\gamma_3} w^{d^3}$ 
is the term in $H^3_{l}$
with the biggest degree with respect to $w$ in $H^3_{l}$;
hence $Q^3$ contains the term $z^{\gamma_3} w^{d^3}$.
Repeating this process,
one can show that $z^{\gamma_n} w^{d^n}$
is the term in $H^n_{l}$
with the biggest degree $d^n$ with respect to $w$ in $H^n_{l}$ inductively;
hence the term $z^{\gamma_n} w^{d^n}$ remains forever.
\end{proof}

The proof of Theorem \ref{main thm of weights for case 3}
also shows that $(\gamma_n, d^n)$ is the vertex of $N(Q^n)$.
Moreover,
we can prove that 
$(\gamma_n, d^n)$ is the vertex of $N(Q^n)$
whose $x$-coordinate is minimum. 
In fact,
we prove that  
the order of $Q^n$ with respect to $z$ is $\gamma_n$,
whereas we showed for Case 2 that 
the order of $Q^n$ with respect to $w$ is $d^n$
in the proof of Theorem \ref{main thm of weights when d>0}.

\begin{proposition} \label{order wrt z}
The order of $Q^n$ with respect to $z$ is $\gamma_n$ for any $n \geq 1$.
\end{proposition}

\begin{proof}
The equality for $n=1$ follows from the setting.

We first show the equality for $n=2$.
Although 
$Q^2(z,w) = \sum b_{ij} \big( p(z) \big)^i \big( q(z,w) \big)^j$,
it is enough to consider the part
$\sum b_{ij} \big( z^{\delta} \big)^{i} \big( z^{\gamma} w^d \big)^{j}$
of $Q^2$ and show that
\[
\delta i + \gamma j \geq \delta \gamma + \gamma d = \gamma_2 
\]
for any $(i,j)$ such that $b_{ij} \neq 0$
since we only interested in the order with respect to $z$.
It follows from Lemma \ref{dyn of R for case 3}
that $i + R(l) j \geq \gamma + R(l) d$,
and so
$\delta i + \gamma j + l dj \geq \delta \gamma + \gamma d + l d^2$ 
by Lemma \ref{R^n}.
Let $(i,j) \neq (\gamma,d)$. 
If $j < d$, then the inequality above implies that
$\delta i + \gamma j > \delta \gamma + \gamma d$
since $ldj < ld^2$.
If $j = d$, then
$\delta i + \gamma j > \delta \gamma + \gamma d$
since $i > \gamma$.
If $j > d$, then
$\delta i + \gamma j > \delta \gamma + \gamma d$
since $i \geq \gamma$.
Therefore, 
$\delta i + \gamma j > \delta \gamma + \gamma d$
for any $(i,j) \neq (\gamma,d)$ and so
the order of $Q^2$ with respect to $z$ is $\gamma_2$.

We next show the equality for $n=3$
by the same strategy.
Although 
$Q^3(z,w) = \sum b_{ij} \big( p^2(z) \big)^i \big( Q^2(z,w) \big)^j$,
it is enough to consider the part
$\sum b_{ij} \big( z^{\delta^2} \big)^{i} \big( z^{\gamma_2} w^{d^2} \big)^{j}$
of $Q^3$ and show that
\[
\delta^2 i + \gamma_2 j \geq \delta^2 \gamma + \gamma_2 d = \gamma_3 
\]
for any $(i,j)$ such that $b_{ij} \neq 0$.
It follows from Lemma \ref{dyn of R for case 3}
that $i + R^2(l) j \geq \gamma + R^2(l) d$,
and so
$\delta^2 i + \gamma_2 j + l d^2 j \geq \delta^2 \gamma + \gamma_2 d + l d^3$ 
by Lemma \ref{R^n}.
If $(i,j) \neq (\gamma,d)$,
then $\delta^2 i + \gamma_2 j > \delta^2 \gamma + \gamma_2 d$
as the same as above.
Therefore, 
the order of $Q^3$ with respect to $z$ is $\gamma_3$.

Repeating this process,
we obtain the required equalities.
\end{proof}

\begin{corollary} \label{cor on N(Q^n) for case 3}
The bidegree
$(\gamma_n, d^n)$ is the vertex of $N(Q^n)$
whose $x$-coordinate is minimum
for any $n \geq 1$. 
\end{corollary}

Therefore,
$N(Q^n)$ is included in the upper-right region
that is surrounded by the line with slope $-l_2^{-1}$
and the vertical line,
which intersect at $(\gamma_n, d^n)$.

Let $- M_n$ be the slope of the line passing through 
$(\gamma_n, d^n)$ and the next vertex if it exist,
let $M_n = \infty$ if $N(Q^n) = D(\gamma_n, d^n)$,
and let $M = M_1$.
Then $M = l_2^{-1}$, and
Theorem \ref{main thm of weights for case 3} and 
Proposition \ref{order wrt z} or Corollary \ref{cor on N(Q^n) for case 3}
imply the following corollary
and Corollary \ref{main cor on attr rates for case 3}.

\begin{corollary}
It follows that $M_n \leq M$ for any $n \geq 1$.
\end{corollary}

\begin{remark} 
If $\gamma > 0$, 
then $\gamma_n \to \infty$
 as $n \to \infty$
since $\delta > 0$ and $d > 0$.
Hence
$c(Q^n) \to \infty$ as $n \to \infty$
if $\gamma > 0$ or $d > 1$.
On the other hand, 
$\min \{ 1, l_2 \} d^n \leq c(Q^n) \leq d^n$
if $\gamma = 0$
and so $c(Q^n) = 1$
if $\gamma = 0$ and $d = 1$.
\end{remark}

\subsection{Next vertices and proof of Theorem \ref{main thm of attr rates for case 3}}

Furthermore,
we can show the existence of the vertex of $N(Q^n)$ 
that is next to $(\gamma_n, d^n)$ in most cases,
and specify it.
Let $(C, D) = (n_{2}, m_{2})$,
\begin{align*}
(C_n, D_n) 
&= ((\delta^{n-1} + \delta^{n-2} d + \cdots + \delta d^{n-2}) \gamma + C d^{n-1}, D d^{n-1}) \\
&= (\gamma_n - (\gamma - C) d^{n-1}, D d^{n-1}) \text{ and} \\
(C_n^*, D_n^*) 
&= ( (\delta^{n-1} + \delta^{n-2} D + \cdots + \delta D^{n-2} + D^{n-1}) C,  D^{n}).
\end{align*}

\begin{proposition} \label{C_n and D_n}
If $\delta > T_{1}$, then
$(C_n, D_n)$ is the vertex of $N(Q^n)$ 
that is next to $(\gamma_n, d^n)$,
$M_n = M$, and
$\delta^n$ is bigger than the $y$-intercept of 
the line passing through $(\gamma_n, d^n)$ and $(C_n, D_n)$
for any $n \geq 1$.
If $\delta = T_{1}$ and $m_2 > 0$, then
$(C_n^*, D_n^*)$ is the vertex of $N(Q^n)$ 
that is next to $(\gamma_n, d^n)$,
$M_n = M$, and
$\delta^n$ coincides with the $y$-intercept of 
the line passing through $(\gamma_n, d^n)$ and $(C_n^*, D_n^*)$
for any $n \geq 1$.
\end{proposition}

\begin{proof}
The proof is similar to that of 
Proposition \ref{A_n and B_n}
for Case 2.
Using the polynomial $H^2_{l_2}$,
we only show how to find out 
$(C_2, D_2)$ and $(C_2^*, D_2^*)$.
If $\delta > T_{1}$,
then it follows from Proposition \ref{H^2 for case 3} that
the desired term is $b_{\gamma d} (a_{\delta} z^{\delta})^{\gamma} (b_{CD} z^C w^D)^d$
and so $(C_2, D_2) = (\delta \gamma + Cd, Dd)$.
If $\delta = T_{1}$,
then it follows from Proposition \ref{H^2 for case 3} that
the desired term is $ b_{CD} (a_{\delta} z^{\delta})^C (b_{CD} z^C w^D)^D$
and so $(C_2^*, D_2^*) = (\delta C + CD,D^2)$.
\end{proof}

Theorem \ref{main thm of weights for case 3} and
this proposition induce  
Theorem \ref{main thm of attr rates for case 3}.

\begin{proof}[Proof of Theorem \ref{main thm of attr rates for case 3}]
We only show the improved parts. \\
(3)
Since $D = m_2 > 0$,
$D_n = D d^{n-1} > 0$ if  $\delta > T_1$, and
$D_n^* = D^n > 0$ if  $\delta = T_1$.
Hence 
$N(Q^n)$ does not contain $(\gamma_n + l_2 d^n, 0)$ and so
$c(Q^n) > \gamma_n + l_2 d^n$.
On the other hand,
$c(Q^n) \leq C_n + D_n < \gamma_n + d^n$ if $\delta > T_1$, and
$c(Q^n) \leq C_n^* + D_n^* < \gamma_n + d^n$ if $\delta = T_1$. \\
(4) 
It follows that 
$D_n = 0$ and $c(Q^n) = C_n = \gamma_n + l_2 d^n$. \\
(5) 
It follows that 
$D_n^* = 0$ and $c(Q^n) = C_n^* = \gamma_n + l_2 d^n$. \\
(6) 
It follows that 
$D_n^* = 0$ and $c(Q^n) > C_n^* = \gamma_n + l_2 d^n$ for any $n \geq n_0$.
\end{proof}

\section{Shape of Newton polygon of $Q^n$ for Case 4}

We finally deal with Case 4 in this section.
Let $s > 2$,
\[
T_k \leq \delta \leq T_{k-1} \text{ for some }2 \leq k \leq s-1,
\] 
\[
(\gamma, d) = (n_k, m_k), \
l_1 = \frac{n_k - n_{k-1}}{m_{k-1} - m_k} \text{ and } 
l_1 + l_2= \frac{n_{k+1} - n_k}{m_k - m_{k+1}}.
\] 
Note that $\delta > d > 0$ and $\gamma > 0$ by the setting.
We first give another interval of weights
which is closely related to those in Section 2.3 and
a summary of the results in Section 6.1.
Some of the results are illustrated 
in terms of the blow-ups in Section 6.2. 
Preparing a lemma in Section 6.3,
we prove 
Theorem \ref{main thm of weights for case 4} below in Section 6.4.
Finally,
we state the claims on the existence of the vertices of $N(Q^n)$
that are previous and next to $(\gamma_n, d^n)$
in most cases
and specify it in Section 6.5,
which induces Theorem \ref{main thm of attr rates for case 4} below.

\subsection{Another interval and Summary of results}

We define the interval $\mathcal{I}_f^{AR}$ as
\[
\mathcal{I}_f^{AR} = 
\left\{ \ l > 0 \ | 
\begin{array}{lcr}
\gamma + ld \leq i + l j
\text{ for any $i$ and $j$ such that } b_{ij} \neq 0
\end{array}
\right\}.
\]
Then
\begin{align*}
\mathcal{I}_f^{AR} &= 
\left\{ \ l > 0 \ | 
\begin{array}{lcr}
\gamma + ld \leq n_j + l m_j
\text{ for any } 1 \leq j \leq s
\end{array}
\right\} \\
&= 
\left\{ \ l > 0 \ \Bigg| 
\begin{array}{lcr}
\gamma + l d \leq n_{k-1} + l m_{k-1} \\ 
\gamma + l d \leq n_{k+1} + l m_{k+1}  
\end{array}
\right\} \\
&=
\left[
\dfrac{\gamma - n_{k-1}}{m_{k-1} - d},
\dfrac{n_{k+1} - \gamma}{d - m_{k+1}} 
\right]
= [l_1, l_1 + l_2].
\end{align*}
As shown in the theorem below,
this interval is suitable for describing the result on $w(Q^n)$
rather than the intervals $\mathcal{I}_f^1$ and $\mathcal{I}_f^2$
and the rectangle $\mathcal{I}_f$ in Section 2.3. 

\begin{theorem} \label{main thm of weights for case 4}
It follows that $Q^n$ contains the term $z^{\gamma_n} w^{d^n}$ and 
$w(Q^n) = \gamma_n + l d^n$
for any $n \geq 1$ and
for any $l \in \mathcal{I}_f^{AR}$,
where $\mathcal{I}_f^{AR} = [l_1, l_1 + l_2]$.
\end{theorem}

Moreover,
$(\gamma_n, d^n)$ is a vertex of $N(Q^n)$ for any $n \geq 1$.
Hence the theorem above implies the following corollary.

\begin{corollary} \label{main cor on attr rates for case 4}
It follows for any $n \geq 1$ that
\begin{enumerate}
\item[$(1)$] $c(Q^n) = \gamma_n + d^n$
if $l_1 \leq 1 \leq l_1 + l_2$, 
\item[$(2)$] $l_1^{-1} \gamma_n + d^n \leq c(Q^n) < \gamma_n + d^n$
if $l_1 > 1$, and
\item[$(3)$] $\gamma_n + (l_1 + l_2) d^n \leq c(Q^n) \leq \gamma_n + d^n$
if $l_1 + l_2 < 1$.
\end{enumerate}
\end{corollary}

Similar to Cases 2 and 3,
by investigating the vertices of $N(Q^n)$
that are previous and next to $(\gamma_n, d^n)$,
we can improve  
the corollary above as follows. 

\begin{theorem} \label{main thm of attr rates for case 4}
It follows for any $n \geq 1$ that
\begin{enumerate}
\item[$(1)$] $c(Q^n) = \gamma_n + d^n$
if $l_1 \leq 1 \leq l_1 + l_2$, 
\item[$(2)$] $l_1^{-1} \gamma_n + d^n \leq c(Q^n) < \gamma_n + d^n$
if $l_1 > 1$, and
\item[$(3)$] $\gamma_n + (l_1 + l_2) d^n \leq c(Q^n) \leq \gamma_n + d^n$
if $l_1 + l_2 < 1$.
\end{enumerate}
More precisely,
it follows for any $n \geq 1$ that
\begin{enumerate}
\item[$(4)$] $l_1^{-1} \gamma_n + d^n < c(Q^n)$ 
             if $l_1 > 1$ and $n_{k-1} > 0$, and
\item[$(5)$] $\gamma_n + (l_1 + l_2) d^n < c(Q^n) < \gamma_n + d^n$ 
             if $l_1 + l_2 < 1$ and $m_{k+1} > 0$. 
\end{enumerate}
Let $l_1 + l_2 < 1$ and $m_{k+1} = 0$.
Then
\begin{enumerate}
\item[$(6)$] $\gamma_n + (l_1 + l_2) d^n = c(Q^n)$ for any $n \geq 1$ if $\delta > T_k$, 
\item[$(7)$] $\gamma_n + (l_1 + l_2) d^n = c(Q^n)$ for any $n \geq 1$ 
        if $\delta = T_k$ and the term $z^{C_n^*} w^{D_n^*}$ remains forever, and
\item[$(8)$] $\gamma_n + (l_1 + l_2) d^n < c(Q^n)$ for any $n \geq n_0$ 
        if $\delta = T_k$ and the term $z^{C_n^*} w^{D_n^*}$ vanishes for some $n_0$.
\end{enumerate}
\end{theorem}

\subsection{Illustration of results in terms of blow-ups}

Assuming that $l_1$ and $l_2^{-1}$ are integers,
we illustrate Theorem \ref{main thm of weights for case 4} and 
Corollary \ref{main cor on attr rates for case 4} in terms of the blow-ups. 
Let $\pi_1 (z,c) = (z, z^{l_1} c)$
and $\tilde{f}_1 = \pi_1^{-1} \circ f \circ \pi_1$
as in Case 2 and
let $\pi_2 (t,w) = (t w^{l_2^{-1}}, w)$
and $\tilde{f}_2 = \pi_2^{-1} \circ \tilde{f}_1 \circ \pi_2$
as in Case 3. 
By Lemma \ref{blow-up lem for case 4},
$N(\tilde{q}_2) = D(\tilde{\gamma}, \tilde{d})$.
Although $\tilde{f}_2$ may not be skew product,
it is close to a monomial map 
by Proposition \ref{blow-up prop for case 4} 
and so $N(\tilde{Q}_2^n)$ has the unique vertex,
that should be $(\tilde{\gamma}_n, \tilde{d}^n)$.
Hence 
$N(\tilde{Q}_2^n) = D(\tilde{\gamma}_n, \tilde{d}^n)$.
This implies that 
$Q^n$ contains the term $z^{\gamma_n} w^{d^n}$, and that
$N(Q^n)$ is included in the upper-right region
that is surrounded by the two lines with slopes $-l_1^{-1}$ and $-(l_1+l_2)^{-1}$,
which intersect at $(\gamma_n, d^n)$. 
Therefore,
we obtain Theorem \ref{main thm of weights for case 4}
and Corollary \ref{main cor on attr rates for case 4} when $l_1$ and $l_2^{-1}$ are integers,
because
$w(Q^n)$ and $c(Q^n)$ are the minimum $x$-intercepts of 
the lines with slopes $-l^{-1}$ and $-1$ 
that intersect $N(Q^n)$. 

\subsection{Preliminary lemma: Affine dynamics on the interval}

Note that $R$ is a contracting function around the fixed point $\alpha$
since $\delta > d$.
Since $\mathcal{I}_f^{AR}$ contains $\alpha$,
we have the following lemma.

\begin{lemma} \label{dyn of R for case 4}
For any $l$ in $\mathcal{I}_f^{AR}$,
it follows that $R^n (l)$ belongs to $\mathcal{I}_f^{AR}$ for any $n \geq 1$ and 
$R^n (l) \to \alpha$ as $n \to \infty$.
More precisely, 
the sequence $\{ R^n (l) \}_{n \geq 1}$ is increasing if $l < \alpha$
and decreasing if $l > \alpha$.
\end{lemma}

\subsection{Proof of Theorem \ref{main thm of weights for case 4}} 

The proof of Theorem \ref{main thm of weights for case 4} is 
almost the same as the proof of Theorem \ref{main thm of weights for case 3}
for Case 3;
we use the polynomial $H^n_l$ 
to show that 
$Q^n$ contains the term $z^{\gamma_n} w^{d^n}$.

\begin{proposition} \label{H^2 for case 4}
Let $l=l_1$ and $(A, B) = (n_{k-1}, m_{k-1})$.
Then
\begin{align*}
H^2_{l_1} (z,w) &= b_{\gamma d} \left( a_{\delta} z^{\delta} \right)^{\gamma} 
           \bigg( \sum_{i+l_1 j= \gamma +l_1 d} b_{ij} z^i w^j \bigg)^d \\
&= b_{\gamma d} \big( a_{\delta} z^{\delta} \big)^{\gamma} 
  \big( b_{\gamma d} z^{\gamma} w^d + \cdots + b_{AB} z^A w^B \big)^d
\text{ if $\delta < T_{k-1}$, and}
\end{align*}
\begin{align*}
H^2_{l_1} (z,w) &= \sum_{I+l_1J= \gamma +l_1 d} b_{IJ}
 \left( a_{\delta} z^{\delta} \right)^{I} \bigg( \sum_{i+l_1 j= \gamma +l_1 d} b_{ij} z^i w^j \bigg)^{J} \\
&= b_{\gamma d} \big( a_{\delta} z^{\delta} \big)^{\gamma} 
  \big( b_{\gamma d} z^{\gamma} w^d + \cdots + b_{AB} z^A w^B \big)^d + \cdots \\
& \ \ \ \ + b_{AB} \big( a_{\delta} z^{\delta} \big)^{A} 
  \big( b_{\gamma d} z^{\gamma} w^d + \cdots + b_{AB} z^A w^B \big)^B
\text{ if $\delta = T_{k-1}$}.
\end{align*}
Let $l = l_1 + l_2$ and $(C, D) = (n_{k+1}, m_{k+1})$.
Then
\begin{align*}
H^2_{l_1 + l_2} (z,w) &=  b_{\gamma d}
 \left( a_{\delta} z^{\delta} \right)^{\gamma} \bigg( \sum_{i+(l_1+l_2) j= \gamma + (l_1+l_2) d} b_{ij} z^i w^j \bigg)^{d} \\
&= b_{\gamma d}\big( a_{\delta} z^{\delta} \big)^{\gamma} 
  \big( b_{\gamma d} z^{\gamma} w^d + \cdots + b_{CD} z^C w^D \big)^d 
\text{ if $\delta > T_{k+1}$, and}
\end{align*}
\begin{align*}
H^2_{l_1 + l_2} (z,w) &= \sum_{I+(l_1+l_2)J= \gamma +(l_1+l_2) d} b_{IJ}
 \left( a_{\delta} z^{\delta} \right)^{I} \bigg( \sum_{i+(l_1+l_2) j= \gamma + (l_1+l_2) d} b_{ij} z^i w^j \bigg)^{J} \\
&= b_{\gamma d} \big( a_{\delta} z^{\delta} \big)^{\gamma} 
  \big( b_{\gamma d} z^{\gamma} w^d + \cdots + b_{CD} z^C w^D \big)^d + \cdots \\
& \ \ \ \ + b_{CD} \big( a_{\delta} z^{\delta} \big)^{C} 
  \big( b_{\gamma d} z^{\gamma} w^d + \cdots + b_{CD} z^C w^D \big)^D
\text{ if $\delta = T_{k+1}$}.
\end{align*}
Let $l_1 < l < l_1 + l_2$.
Then
\[
H^2_{l} (z,w) = b_{\gamma d} \big( a_{\delta} z^{\delta} \big)^{\gamma} \big( b_{\gamma d} z^{\gamma} w^d \big)^d.
\]
\end{proposition}

\begin{proof}[Proof of Theorem \ref{main thm of weights for case 4}] 
Fix any $l$ in $\mathcal{I}_f^{AR}$.
The inequality $w(Q^n) \geq \gamma_n + l d^n$ follows from
Lemma \ref{dyn of R for case 4} and
the same arguments as the proof of 
Theorem \ref{main thm of weights when d>0} for Case 2.


It follows from Proposition \ref{H^2 for case 4} that
$z^{\gamma_2} w^{d^2}$ is the term in $H^2_l$ with 
the smallest and/or biggest degree with respect to $w$ in $H^2_{l}$
and so $Q^2$ contains the term $z^{\gamma_2} w^{d^2}$.
More precisely,
$d^2$ is the smallest degree with respect to $w$ in $H^2_{l}$ if $l = l_1$,
the biggest degree with respect to $w$ in $H^2_l$ if $l = l_1 + l_2$,
and $H^2_{l}$ consists only of the term $z^{\gamma_2} w^{d^2}$
if $l_1 < l < l_1 + l_2$.
Applying the same argument inductively,
one can show for any $n \geq 3$
that $z^{\gamma_n} w^{d^n}$ is the term in $H^n_l$ with
the smallest and/or biggest degree with respect to $w$ in $H^n_l$
and so 
$Q^n$ contains the term $z^{\gamma_n} w^{d^n}$.
\end{proof}

The proof above implies the following.

\begin{corollary} \label{vertex for case 4}
The bidegree
$(\gamma_n, d^n)$ is a vertex of $N(Q^n)$ for any $n \geq 1$.
\end{corollary}

Therefore,
$N(Q^n)$ is included in the upper-right region
that is surrounded by the two lines with slopes $-l_1^{-1}$ and $-(l_1+l_2)^{-1}$,
which intersect at $(\gamma_n, d^n)$. 

Let $- M_n(l_1)$ and $- M_n(l_1 + l_2)$  be the slopes of 
the line passing through 
$(\gamma_n, d^n)$ and the previous vertex
and the line passing through 
$(\gamma_n, d^n)$ and the next vertex,
respectively.
Let $M(l_1) = M_1(l_1)$ and $M(l_1 + l_2) = M_1(l_1 + l_2)$.
Then $M(l_1) = l_1^{-1}$ and $M(l_1 + l_2) = (l_1 + l_2)^{-1}$.
Theorem \ref{main thm of weights for case 4}
and Corollary \ref{vertex for case 4}
imply the following corollary and 
Corollary \ref{main cor on attr rates for case 4}.

\begin{corollary}
It follows that 
$M_n (l_1) \geq M(l_1)$ and
$M_n (l_1 + l_2) \leq M(l_1 + l_2)$
for any $n \geq 1$.
\end{corollary}

\subsection{Claims on previous and next vertices}

As in Cases 2 and 3,
we can show the existence of the vertices of $N(Q^n)$ 
that are previous and next to $(\gamma_n, d^n)$ in most cases
and specify them.
Let $(A, B) = (n_{k-1}, m_{k-1})$ and $(C, D) = (n_{k+1}, m_{k+1})$.
Let $(A_n, B_n)$ and $(A_n^*, B_n^*)$ be the same as Case 2 and
let $(C_n, D_n)$ and $(C_n^*, D_n^*)$ be the same as Case 3.

\begin{proposition}
If $\delta < T_{k-1}$, then
$(A_n, B_n)$ is the vertex of $N(Q^n)$ 
that is previous to $(\gamma_n, d^n)$,
$M_n(l_1) = M(l_1)$, and
$\delta^n$ is smaller than the $y$-intercept of 
the line passing through $(\gamma_n, d^n)$ and $(A_n, B_n)$ for any $n \geq 1$.
If $\delta = T_{k-1}$, then
$(A_n^*, B_n^*)$ is the vertex of $N(Q^n)$ 
that is previous to $(\gamma_n, d^n)$,
$M_n(l_1) = M(l_1)$, and
$\delta^n$ coincides with the $y$-intercept of 
the line passing through $(\gamma_n, d^n)$ and $(A_n^*, B_n^*)$ for any $n \geq 1$.
\end{proposition}

\begin{proposition}
If $\delta > T_{k}$, then
$(C_n, D_n)$ is the vertex of $N(Q^n)$ 
that is next to $(\gamma_n, d^n)$,
$M_n(l_1 + l_2) = M(l_1 + l_2)$, and 
$\delta^n$ is bigger than the $y$-intercept of 
the line passing through $(\gamma_n, d^n)$ and $(C_n, D_n)$ for any $n \geq 1$.
If $\delta = T_{k}$ and $m_{k+1} > 0$, then
$(C_n^*, D_n^*)$ is the vertex of $N(Q^n)$ 
that is next to $(\gamma_n, d^n)$,
$M_n(l_1 + l_2) = M(l_1 + l_2)$, and 
$\delta^n$ coincides with the $y$-intercept of 
the line passing through $(\gamma_n, d^n)$ and $(C_n^*, D_n^*)$ for any $n \geq 1$.
\end{proposition}

The proofs of these propositions are the same as those of
Propositions \ref{A_n and B_n} and \ref{C_n and D_n}
for Cases 2 and 3, and
these propositions induce 
Theorem \ref{main thm of attr rates for case 4}
as in Cases 2 and 3.

\section{Estimates on attraction rate of $f^n$} 

In this last section
we derive inequalities on $c(f^n)$
from the inequalities on $c(Q^n)$ in Theorem \ref{main thm on attr rates}.
We first give a summary of our results on the inequalities on $c(f^n)$ in Section 7.1.
Detailed estimates and explanations 
for Case 1, Case 2 when $d>0$, Case 2 when $d=0$, Case 3 and Case 4
are given in Sections 7.2, 7.3, 7.4, 7.5 and 7.6,
respectively.

\subsection{Summary on inequalities on $c(f^n)$}

It is clear that $c(f^n) \leq \delta^n$ since
$c(f^n) = \min \{ c(p^n), c(Q^n) \} = \min \{ \delta^n, c(Q^n) \} \leq \delta^n$.
Theorem \ref{main thm on attr rates} induces
the following inequalities on $c(f^n)$.

\begin{theorem} 
Let $\gamma d > 0$.
Then it follows for any $n \geq 1$ that
\begin{enumerate}
\item[$(1)$] $\alpha \delta^n \leq c(f^n) < \delta^n$ if $\delta > d$ and $\alpha < 1$, and 
\item[$(2)$] $c(f^n) = \delta^n$ if $\delta > d$ and $\alpha \geq 1$ or if $\delta \leq d$. 
\end{enumerate}
On the other hand,
it follows for any $n \geq 1$ that
\begin{enumerate}
\item[$(3)$] $c(f^n) = \min \{ 1, \gamma/\delta \} \delta^{n}$ if $d = 0$ and $s = 1$,
\item[$(4)$] $c(f^n) = \min \{ \delta^n, d^n \}$ if $\gamma = 0$ and $s = 1$,
\item[$(5)$] $\min \{ 1, \gamma/\delta, l_1^{-1} \gamma/\delta \} \delta^{n} \leq c(f^n) \leq \delta^n$ 
             if $d = 0$ and $s > 1$, and
\item[$(6)$] $\min \{ 1, l_2 \} d^n \leq c(f^n) \leq d^n$ if $\gamma = 0$ and $s > 1$.
\end{enumerate}
In particular,
$c_{\infty} = \delta$ if $\gamma > 0$, and
$c_{\infty} = \min \{ \delta, d \}$ if $\gamma = 0$.
\end{theorem} 

This theorem implies 
Corollary \ref{main cor on attr rates}.
Note that the cases (3) and (4) occur only for Case 1,
the case (5) occurs only for Case 2, and 
the case (6) occurs only for Case 3.
We give improved versions of (5) in Section 7.4
by applying Theorems \ref{main thm of attr rates when d=0 and delta < T}
and \ref{main thm of attr rates when d=0 and delta = T}.

\subsection{Estimates on $c(f^n)$ for Case 1}

For Case 1, 
we may assume that $f$ is a monomial map
in order to estimate the attraction rates.
Let $f(z,w) = (z^{\delta}, z^{\gamma} w^d)$,
where $\delta \geq 1$, $\gamma \geq 0$, $d \geq 0$ and $\gamma + d \geq 1$.
Then $f^n(z,w) = (z^{\delta^n}, z^{\gamma_n} w^{d^n})$ and so
$c(p^n) = \delta^n$ and $c(Q^n) = \gamma_n + d^n$,
where
\[
\gamma_n 
= \gamma \dfrac{\delta^n - d^n}{\delta - d}
= \alpha (\delta^n - d^n)
\text{ if } \delta \neq d,
\text{ and }
\gamma_n
= n \gamma \delta^{n-1}
\text{ if } \delta = d
\]
since $\gamma_n = \gamma (\delta^{n-1} + \delta^{n-2} d + \delta^{n-3} d^2 + \cdots + d^{n-1})$.

First, 
let $\gamma d > 0$.
If $\delta > d$, 
then $\gamma_n + d^n = \alpha \delta^n + (1 - \alpha)d^n \sim \alpha \delta^n$,
where $\alpha > 0$, and 
\begin{equation*}
\begin{cases}
\ \alpha \delta^n < \gamma_n + d^n < \alpha \delta^n + (1 - \alpha) \delta^n = \delta^n
 & \text{ if } \alpha < 1, \\
\ \delta^n = \alpha \delta^n + (1 - \alpha) \delta^n < \gamma_n + d^n < \alpha \delta^n
 & \text{ if } \alpha > 1, \\
\ \gamma_n + d^n = \delta^n 
 & \text{ if } \alpha = 1.
\end{cases}
\end{equation*}
Here the notation $A_n \sim B_n$ means that
the ratio of $A_n$ and $B_n$ tends to $1$ as $n \to \infty$.
Since $c(f^n) = \min \{ \delta^n, c(Q^n) \}$,
we can summarize the estimates on $c(Q^n)$ and $c(f^n)$ 
as follows.

\vspace{-4mm}
\begin{table}[H]
\caption{Estimates for Case 1 when $\gamma d>0$ and $\delta > d$}
\begin{center}
{\begin{tabular}{|c|c|c|} \hline
\rule[-5pt]{0pt}{18pt} $\alpha < 1$ & $\alpha \delta^n < c(Q^n) < \delta^n$ & $\alpha \delta^n < c(f^n) < \delta^n$ \\ \hline
\rule[-5pt]{0pt}{18pt} $\alpha > 1$ & $\delta^n < c(Q^n) < \alpha \delta^n$ & $c(f^n) = \delta^n$ \\ \hline
\rule[-5pt]{0pt}{18pt} $\alpha = 1$ & $c(Q^n) = \delta^n$ & $c(f^n) = \delta^n$ \\ \hline
\end{tabular}}
\end{center}
\end{table}

\noindent
If $\delta < d$,
then $\gamma_n + d^n = (- \alpha + 1) d^n - (- \alpha) \delta^n \sim (- \alpha + 1) d^n$,
where $\alpha < 0$, and 
$d^n < \gamma_n + d^n < (- \alpha + 1) d^n$.
If $\delta = d$, 
then $\gamma_n + d^n = n \gamma \delta^{n-1} + \delta^n \sim n \gamma \delta^{n-1}$
and $\delta^n = d^n < \gamma_n + d^n$.
Consequently,
we can summarize the estimates on $c(Q^n)$
for Case 1 when $\gamma d>0$ as follows.

\vspace{-4mm}
\begin{table}[H]
\caption{Estimates on $c(Q^n)$ for Case 1 when $\gamma d>0$}
\begin{center}
{\begin{tabular}{|c|c|c|} \hline
\rule[-5pt]{0pt}{18pt} $\delta > d$ &  
                           $c(Q^n) \sim \alpha \delta^n$ & 
                           $\min \{ \alpha, 1 \} \delta^n \leq c(Q^n) \leq \max \{ \alpha, 1 \} \delta^n$ \\ \hline
\rule[-5pt]{0pt}{18pt} $\delta < d$ &   
                           $c(Q^n) \sim (- \alpha + 1) d^n$ & 
                           $d^n < c(Q^n) < (- \alpha + 1) d^n$ \\ \hline
\rule[-5pt]{0pt}{18pt} $\delta = d$ &   
                           $c(Q^n) \sim n \gamma \delta^{n-1}$ & 
                           $\delta^n < c(Q^n) = n \gamma \delta^{n-1} + \delta^n$ \\ \hline
\end{tabular}}
\end{center}
\end{table}

\noindent
Since $c(f^n) = \min \{ \delta^n, c(Q^n) \}$,
we obtain the following estimates on $c(f^n)$.

\vspace{-4mm}
\begin{table}[H]
\caption{Estimates on $c(f^n)$ for Case 1 when $\gamma d>0$}
\begin{center}
{\begin{tabular}{|c|c|} \hline
\rule[-5pt]{0pt}{18pt} $\delta > d$ &                            
                           $\min \{ \alpha, 1 \} \delta^n \leq c(f^n) \leq \delta^n$ \\ \hline
\rule[-5pt]{0pt}{18pt} $\delta < d$ &   
                           $c(f^n) = \delta^n$ \\ \hline
\rule[-5pt]{0pt}{18pt} $\delta = d$ &   
                           $c(f^n) = \delta^n$ \\ \hline
\end{tabular}}
\end{center}
\end{table}

Next, 
let $\gamma d = 0$.
If $\gamma = 0$,
then $c(Q^n) = d^n$ and so $c(f^n) = \min \{ \delta^n, d^n \}$.
If $d = 0$,
then $c(Q^n) = \gamma_n = \gamma \delta^{n-1}$ and so $c(f^n) = \min \{ \delta, \gamma \} \delta^{n-1}$.

\vspace{-4mm}
\begin{table}[H]
\caption{Estimates for Case 1 when $\gamma d = 0$}
\begin{center}
{\begin{tabular}{|c|c|c|} \hline
\rule[-5pt]{0pt}{18pt} $\gamma = 0$ & $c(Q^n) = d^n$ & 
                           $c(f^n) = \min \{ \delta^n, d^n \}$ \\ \hline
\rule[-5pt]{0pt}{18pt} $d = 0$ & $c(Q^n) = \gamma \delta^{n-1}$ & 
                           $c(f^n) = \min \{ \delta, \gamma \} \delta^{n-1}$ \\ \hline
\end{tabular}}
\end{center}
\end{table}

\subsection{Estimates on $c(f^n)$ for Case 2 when $d>0$}

Recall that $\gamma > 0$ by the setting.
If $d > 0$ and $l_1 \leq 1$,
then $c(Q^n) = \gamma_n + d^n$ 
by Theorem \ref{main thm on attr rates}
and so
we have the same estimates as in Case 1. 
Let $d > 0$ and $l_1 > 1$.
Then $l_1^{-1} \gamma_n + d^n \leq c(Q^n) \leq \gamma_n + d^n$
by Theorem \ref{main thm on attr rates} and 
$\l_1^{-1} \gamma_n + d^n
= l_1^{-1} \alpha \delta^n + \left( 1 - l_1^{-1} \alpha \right) d^n$.
If $\delta > d$,
then $\mathcal{I}_f = [ l_1, \alpha ]$ and $\alpha \geq l_1 > 1$.
Hence $1 - l_1^{-1} \alpha \leq 0$ and so
$\delta^n
= l_1^{-1} \alpha \delta^n + \left( 1 - l_1^{-1} \alpha \right) \delta^n
\leq l_1^{-1} \gamma_n + d^n
\leq c(Q^n) \leq \gamma_n + d^n < \alpha \delta^n$.
If $\delta < d$,
then
$d^n < l_1^{-1} \gamma_n + d^n
\leq c(Q^n) \leq \gamma_n + d^n < (- \alpha + 1) d^n$.
If $\delta = d$,
then
$\delta^n = d^n < l_1^{-1} \gamma_n + d^n
\leq c(Q^n) \leq \gamma_n + d^n = n \gamma \delta^{n-1} + \delta^n$.
Therefore,
we can summarize the estimates on $c(Q^n)$ and $c(f^n)$ 
for Case 2 when $d>0$ and $l_1 > 1$ as follows,
which are almost the same as Case 1.  

\vspace{-4mm}
\begin{table}[H]
\caption{Estimates for Case 2 when $d>0$ and $l_1 > 1$}
\begin{center}
{\begin{tabular}{|c|c|c|} \hline
\rule[-5pt]{0pt}{18pt} $\delta > d$ (and $\alpha > 1$) & 
                           $\delta^n \leq c(Q^n) < \alpha \delta^n$ & 
                           $c(f^n) = \delta^n$ \\ \hline
\rule[-5pt]{0pt}{18pt} $\delta < d$ &   
                           $d^n < c(Q^n) < (- \alpha + 1) d^n$ & 
                           $c(f^n) = \delta^n$ \\ \hline
\rule[-5pt]{0pt}{18pt} $\delta = d$ &   
                           $\delta^n < c(Q^n) \leq n \gamma \delta^{n-1} + \delta^n$ & 
                           $c(f^n) = \delta^n$ \\ \hline
\end{tabular}}
\end{center}
\end{table}

\noindent
In particular, 
$\delta^n \leq c(Q^n)$ and $c(f^n) = \delta^n$ if $d > 0$ and $l_1 > 1$.

\subsection{Estimates on $c(f^n)$ for Case 2 when $d=0$}

Let $d = 0$. 
If $\delta < T_{s-1}$,
then
\begin{equation*}
\begin{cases}
\ c(Q^n) = \gamma_n
 & \text{ if } l_1 \leq 1, \\
\ l_1^{-1} \gamma_n \leq c(Q^n) \leq \gamma_n
 & \text{ if } l_1 > 1
\end{cases}
\end{equation*}
by Theorem \ref{main thm of attr rates when d=0 and delta < T}.
If $\delta = T_{s-1}$,
then
\begin{equation*}
\begin{cases}
\ \gamma_n \leq c(Q^n) \leq l_1^{-1} \gamma_n
 & \text{ if } l_1 \leq 1, \\
\ l_1^{-1} \gamma_n \leq c(Q^n) < \gamma_n
 & \text{ if } l_1 > 1 
\end{cases}
\end{equation*}
by Theorem \ref{main thm of attr rates when d=0 and delta = T}.
Since $\gamma_n = \gamma \delta^{n-1}$ if $d = 0$, 
we can summarize the estimates on $c(Q^n)$ and $c(f^n)$ 
for Case 2 when $d=0$ as follows.

\vspace{-4mm}
\begin{table}[H]
{\small
\caption{Estimates for Case 2 when $d=0$ and $\delta < T_{s-1}$}
\begin{center}
{\begin{tabular}{|c|c|c|} \hline
\rule[-5pt]{0pt}{18pt} $\delta < T_{s-1}$ & $c(Q^n)$ & $c(f^n)$ \\ \hline
\rule[-5pt]{0pt}{18pt} $l_1 \leq 1$ & $c(Q^n) = \gamma \delta^{n-1}$ & 
                           $c(f^n) = \min \{ \delta, \gamma \} \delta^{n-1}$ \\ \hline
\rule[-5pt]{0pt}{18pt} $l_1 > 1$ & $l_1^{-1} \gamma \delta^{n-1} \leq c(Q^n) \leq \gamma \delta^{n-1}$ & 
                           $\min \{ \delta, l_1^{-1} \gamma \} \delta^{n-1} \leq c(f^n) \leq \min \{ \delta, \gamma \} \delta^{n-1}$ \\ \hline
\end{tabular}}
\end{center}
}
\end{table}

\vspace{-8mm}
\begin{table}[H]
{\small
\caption{Estimates for Case 2 when $d=0$ and $\delta = T_{s-1}$}
\begin{center}
{\begin{tabular}{|c|c|c|} \hline
\rule[-5pt]{0pt}{18pt} $\delta = T_{s-1}$ & $c(Q^n)$ & $c(f^n)$ \\ \hline
\rule[-5pt]{0pt}{18pt} $l_1 \leq 1$ & $\gamma \delta^{n-1} \leq c(Q^n) \leq l_1^{-1} \gamma \delta^{n-1}$ & 
                           $\min \{ \delta, \gamma \} \delta^{n-1} \leq c(f^n) \leq \min \{ \delta, l_1^{-1} \gamma \} \delta^{n-1}$ \\ \hline
\rule[-5pt]{0pt}{18pt} $l_1 > 1$ & $l_1^{-1} \gamma \delta^{n-1} \leq c(Q^n) < \gamma \delta^{n-1}$ & 
                           $\min \{ \delta, l_1^{-1} \gamma \} \delta^{n-1} \leq c(f^n) \leq \min \{ \delta, \gamma \} \delta^{n-1}$ \\ \hline
\end{tabular}}
\end{center}
}
\end{table}

\noindent
In particular, 
$\min \{ 1, \gamma/\delta, l_1^{-1} \gamma/\delta \} \delta^{n} \leq c(f^n) \leq \delta^n$ 
if $d = 0$.
We remark that this rough inequality also follows from 
Theorem \ref{main thm on attr rates}.

\subsection{Estimates on $c(f^n)$ for Case 3}

Recall that
$\delta > d > 0$ if $\gamma > 0$, and
$\delta \geq d > 0$ if $\gamma = 0$
by the setting.
If $l_2 \geq 1$,
then $c(Q^n) = \gamma_n + d^n$ 
by Theorem \ref{main thm on attr rates}
and so
we have the same estimates as in Case 1. 
Let $l_2 < 1$.
Then $\gamma_n + l_2 d^n \leq c(Q^n) \leq \gamma_n + d^n$
by Theorem \ref{main thm on attr rates}.
If $\gamma > 0$,
then $\mathcal{I}_f = [ \alpha, l_2 ]$ and $0 < \alpha \leq l_2 < 1$.
Hence 
$\alpha \delta^n \leq \alpha \delta^n + (l_2 - \alpha) d^n 
= \gamma_n + l_2 d^n \leq c(Q^n) \leq \gamma_n + d^n < \delta^n$. 
If $\gamma = 0$,
then $l_2 d^n \leq c(Q^n) \leq d^n$. 
Therefore,
we can summarize the estimates on $c(Q^n)$ and $c(f^n)$ 
for Case 3 when $l_2 < 1$ as follows.

\vspace{-4mm}
\begin{table}[H]
\caption{Estimates for Case 3 when $l_2 < 1$}
\begin{center}
{\begin{tabular}{|c|c|c|} \hline
\rule[-5pt]{0pt}{18pt} $\gamma > 0$ (and $\delta > d$) & 
                           $\alpha \delta^n \leq c(Q^n) < \delta^n$ & 
                           $\alpha \delta^n \leq c(f^n) < \delta^n$ \\ \hline
\rule[-5pt]{0pt}{18pt} $\gamma = 0$  (and $\delta \geq d$) & 
                           $l_2 d^n \leq c(Q^n) \leq d^n$ & 
                           $l_2 d^n \leq c(f^n) \leq d^n$ \\ \hline
\end{tabular}}
\end{center}
\end{table}

\subsection{Estimates on $c(f^n)$ for Case 4}

Recall that
$\delta > d > 0$ and $\gamma > 0$
by the setting
and that $l_1 \leq \alpha \leq l_1 + l_2$.
If $l_1 \leq 1 \leq l_1 + l_2$,
then $c(Q^n) = \gamma_n + d^n$ 
by Theorem \ref{main thm on attr rates}
and so
we have the same estimates as Case 1. 
If $l_1 > 1$,
then $l_1^{-1} \gamma_n + d^n \leq c(Q^n) \leq \gamma_n + d^n$
by Theorem \ref{main thm on attr rates}
and so $\delta^n \leq c(Q^n) < \alpha \delta^n$
since $\alpha \geq l_1 > 1$.
If $l_1 + l_2 < 1$,
then $\gamma_n + (l_1 + l_2) d^n \leq c(Q^n) \leq \gamma_n + d^n$
by Theorem \ref{main thm on attr rates}
and so
$\alpha \delta^n \leq \alpha \delta^n + (l_1 + l_2 - \alpha) d^n 
= \gamma_n + (l_1 + l_2) d^n \leq c(Q^n) \leq \gamma_n + d^n
< \delta^n$
since $\alpha \leq l_1 + l_2 < 1$. 
We can classify the estimates on $c(Q^n)$ and $c(f^n)$ 
for Case 4 in terms of $\alpha$ 
and obtain the following summary,
which is almost the same as Case 1 
when $\delta > d > 0$ and $\gamma > 0$.

\vspace{-4mm}
\begin{table}[H]
\caption{Estimates for Case 4}
\begin{center}
{\begin{tabular}{|c|c|c|} \hline
\rule[-5pt]{0pt}{18pt} $\alpha < 1$ & $\alpha \delta^n \leq c(Q^n) < \delta^n$ & 
                           $\alpha \delta^n \leq c(f^n) < \delta^n$ \\ \hline
\rule[-5pt]{0pt}{18pt} $\alpha > 1$ & $\delta^n \leq c(Q^n) < \alpha \delta^n$ & $c(f^n) = \delta^n$ \\ \hline
\rule[-5pt]{0pt}{18pt} $\alpha = 1$ & $c(Q^n) = \delta^n$ & $c(f^n) = \delta^n$ \\ \hline
\end{tabular}}
\end{center}
\end{table}



\begin{thebibliography}{9}

\bibitem{ab}
  \textsc{M. Astorg and L. Boc-Thaler}, 
  \textit{Dynamics of skew-products tangent to the identity},
  preprint arXiv:2204.02644 (2022).

\bibitem{abp}
  \textsc{M. Astorg, L. Boc-Thaler and H. Peters}, 
  \textit{Wandering domains arising from Lavaurs maps with Siegel disks},
  Analysis \& PDE \textbf{16} (2023), 35-88.
\bibitem{abdpr}
  \textsc{M. Astorg, X. Buff, R. Dujardin, H. Peters and J. Raissy}, 
  \textit{A two-dimensional polynomial mapping with a wandering Fatou component},
  Annals of Mathematics \textbf{184} (2016), 263-313.

\bibitem{d}
  \textsc{R. Dujardin}, 
  \textit{Non-density of stability for holomorphic mappings on $\mathbb{P}^{k}$},
  J. \`{E}c. polytech. Math. \textbf{4} (2017), 813-843.

\bibitem{fg}
  \textsc{C. Favre and V. Guedj}, 
  \textit{Dynamique des applications rationnelles des espaces multiprojectifs},
  Indiana Univ. Math. J. \textbf{50} (2001), 881-934. 

\bibitem{fj}
  \textsc{C. Favre and M. Jonsson}, 
  \textit{Eigenvaluations},
  Ann. Sci. \'Ecole Norm. Sup. \textbf{40} (2007), pp. 309-349.

\bibitem{gr}
  \textsc{W. Gignac and M. Ruggiero}, 
  \textit{Growth of attraction rates for iterates of a superattracting germ in dimension two},
  Indiana Univ. Math. J. \textbf{63} (2014), pp. 1195-1234.

\bibitem{j}
  \textsc{M. Jonsson}, 
  \textit{Dynamics of polynomial skew products on $\mathbf{C}^{2}$},
  Math. Ann. \textbf{314} (1999), 403-447. 

\bibitem{l}
  \textsc{K. Lilov}, 
  \textit{Fatou theory in two dimensions},
  PhD thesis, University of Michigan, 2004.

\bibitem{pr}
  \textsc{H. Peters and J. Raissy}, 
  \textit{Fatou components of elliptic polynomial skew products},
  Ergodic Theory Dynam. Systems \textbf{39} (2019), 2235-2247. 

\bibitem{ps}
  \textsc{H. Peters and I. M. Smit}, 
  \textit{Fatou components of attracting skew-products},
  J. Geom. Anal. \textbf{28} (2018), 84-110.

\bibitem{t}
  \textsc{J. Taflin}, 
  \textit{Blenders near polynomial product maps of $\mathbb{C}^{2}$},
  J. Eur. Math. Soc. \textbf{23} (2021), 3555-3589. 

\bibitem{u}
  \textsc{K. Ueno}, 
  \textit{A construction of B\"{o}ttcher coordinates for holomorphic skew products},
  Nonlinearity \textbf{32} (2019), 2694-2720. 

\bibitem{u - psh funs}
  \textsc{---}, 
  \textit{Dynamics of superattracting skew products on the attracting basins: 
  B\"{o}ttcher coordinates and plurisubharmonic functions},
  preprint arXiv:2304.09457 (2023).


\end{thebibliography}
\end{document}